\documentclass[11pt,reqno]{article}
\usepackage{amssymb}
\usepackage{amsmath}
\usepackage{amsthm}
\newcommand{\Xbar}{$\Bar{\bf X}$}
\newcommand{\mXbar}{\mbox{\Xbar}}
\newcommand{\Ybar}{$\bar{\bf Y}$}
\newcommand{\mYbar}{\mbox{\Ybar}}

\newcommand{\bSigma}{\boldmath$\Sigma$}
\newcommand{\mbSigma}{\mbox{\bSigma}}
\newcommand{\bLambda}{\boldmath$\Lambda$}
\newcommand{\mbLambda}{\mbox{\bLambda}}
\newcommand{\bGamma}{\boldmath$\Gamma$}
\newcommand{\mbGamma}{\mbox{\bGamma}}

\newcommand{\btheta}{\boldmath$\theta$}
\newcommand{\mbtheta}{\mbox{\btheta}}
\newcommand{\bTheta}{\boldmath$\Theta$}
\newcommand{\mbTheta}{\mbox{\bTheta}}

\newcommand{\bPsi}{\boldmath$\Psi$}
\newcommand{\mbPsi}{\mbox{\bPsi}}
\newcommand{\bUpsilon}{\boldmath$\Upsilon$}
\newcommand{\mbUpsilon}{\mbox{\bUpsilon}}

\newcommand{\bmu}{\boldmath$\mu$}
\newcommand{\mbmu}{\mbox{\bmu}}
\newcommand{\bnu}{\boldmath$\nu$}
\newcommand{\mbnu}{\mbox{\bnu}}

\newcommand{\mbM}{{\bf M}}
\newcommand{\mbN}{{\bf N}}
\newcommand{\mbS}{{\bf S}}
\newcommand{\mbT}{{\bf T}}
\newcommand{\mbZ}{{\bf Z}}
\newcommand{\mbg}{\boldsymbol{g}}
\newcommand{\mbw}{\boldsymbol{w}}
\DeclareMathOperator{\tr}{tr}

\begin{document}

\thispagestyle{empty}
\newtheorem{theorem}{\indent {\bf Theorem}}
\newtheorem{prop}{\indent {\bf Proposition}}
\newtheorem{lemma}{\indent {\bf Lemma}}
\renewcommand{\proofname}{\hspace*{\parindent}{\bf Proof.}}

\begin{center}
 {\bf Admissibility of invariant tests for means with covariates}

\vspace{0.3cm}
\indent Ming-Tien Tsai

\vspace{0.3cm}
\noindent {\it Institute of Statistical Science, Academia Sinica, Taipei 11529, R.O.C.}
\end{center}

\vspace{0.5cm}
\begin{center}
\small
\parbox{13.5cm}{\sloppy \ \ \ \,\,\, For a multinormal distribution with a $p$-dimensional mean
vector ${\mbtheta}$ and an arbitrary unknown dispersion matrix ${\mbSigma}$, Rao ([9], [10])
proposed two tests for the problem of testing $ H_{0}:{\mbtheta}_{1} = {\bf 0}, {\mbtheta}_{2} = {\bf 0},
{\mbSigma}~ \hbox{unspecified},~\hbox{versus}~H_{1}:{\mbtheta}_{1} \ne {\bf 0}, {\mbtheta}_{2}
={\bf 0}, {\mbSigma}~\hbox{unspecified}$, where
${\mbtheta}^{'}=({\mbtheta}^{'}_{1},{\mbtheta}^{'}_{2})$. These tests are referred to as Rao's
$W$-test (likelihood ratio test) and Rao's $U$-test (union-intersection test), respectively.
This work is inspired by the well-known work of Marden and Perlman [6] who claimed that
Hotelling's $T^{2}$-test is admissible while Rao's $U$-test is inadmissible. Both Rao's $U$-test
and Hotelling's $T^{2}$-test can be constructed by applying the union-intersection principle that
incorporates the information ${\mbtheta}_{2}={\bf 0}$ for Rao's $U$-test statistic but does not
incorporate it for Hotelling's $T^{2}$-test statistic. Rao's $U$-test is believed to exhibit some
optimal properties. Rao's $U$-test is shown to be admissible by fully incorporating the information
${\mbtheta}_{2}={\bf 0}$, but Hotelling's $T^{2}$-test is inadmissible.
}
\end{center}

\vspace{0.3cm}
\begin{center}
\small
\parbox{14cm} {\small{\it AMS subject classification}. 62C15, 62H15.}
\end{center}

\begin{center}
\small
\parbox{14cm} {\small
 }
\end{center}

\vspace{0.5cm}
\indent Abbreviated title: Admissibility of invariant tests for means with covariates.
\pagebreak

\vspace{0.3cm}
\def \theequation{1.\arabic{equation}}
\setcounter{equation}{0}
\noindent {\bf 1. Introduction}
\vspace{0.3cm}

\indent Let \{${\bf X}_{i}; 1 \leq i \leq n$\} be independent and identically distributed
random vectors (i.i.d.r.v.) with a $p$-variate normal distribution with mean vector
${\mbtheta}$ and dispersion matrix ${\mbSigma}$, where ${\mbSigma}$ is assumed to be
positive definite (p.d.). Partition ${\mbtheta}$ and ${\mbSigma}$ as
\begin{eqnarray}
  {\mbtheta}=\left(
  \begin{array}{c}
    {\mbtheta}_{1}\\
    {\mbtheta}_{2}
  \end{array}
  \right) ~~\hbox{and}~~
  {\mbSigma} = \left[
  \begin{array}{cc}
    {\mbSigma}_{11}& {\mbSigma}_{12}\\
    {\mbSigma}_{21}& {\mbSigma}_{22}
  \end{array}
  \right],
\end{eqnarray}
where ${\mbtheta}_{1}:p_{1}\times 1,{\mbtheta}_{2}:p_{2}\times 1,{\mbSigma}_{11}:p_{1}
\times p_{1}, {\mbSigma}_{22}:p_{2}\times p_{2}, p_{1}+p_{2}=p$, $0<p_{2}<p.$
The problem of interest is to test
\begin{eqnarray} \nonumber
 & &H_{0}: {\mbtheta}_{1} = {\bf 0},
 {\mbtheta}_{2} = {\bf 0}, ~{\mbSigma}~ \hbox{unspecified} \\
 \hbox{versus} & & \\
  & & H_{1}: {\mbtheta}_{1} \ne {\bf 0},
  {\mbtheta}_{2}={\bf 0}, ~{\mbSigma} ~\hbox{unspecified}. \nonumber
\end{eqnarray}

\indent For every $n$ $(\geq 2)$, let
\begin{eqnarray}
  {\mXbar} = n^{-1} \sum_{i=1}^{n} {\bf X}_{i}~~\hbox{and}~~{\bf S}=
  \sum_{i=1}^{n} ({\bf X}_{i} - {\mXbar}) ({\bf X}_{i}-{\mXbar})',
\end{eqnarray}
and express Hotelling's $T^{2}$-statistic as
\begin{eqnarray}
T^{2}=n(n-1){\mXbar}^{'}{\bf S}^{-1}{\mXbar}.
\end{eqnarray}

\indent Partition ${\mXbar}$ and ${\bf S}$ similarly as in (1.1), and define
\begin{eqnarray}
  {\mXbar}_{1:2} = {\mXbar}_{1}-{\bf S}_{12}{\bf S}_{22}^{-1}{\mXbar}_{2}, \\
  {\bf S}_{11:2}={\bf S}_{11}-{\bf S}_{12}{\bf S}_{22}^{-1}{\bf S}_{21}.
\end{eqnarray}
For the problem (1.2), Rao ([9], [10]) proposed two test statistics which are of the forms
\begin{eqnarray}
  W ={\frac{n(n-1){\mXbar}'_{1:2}{\bf S}^{-1}_{11:2}{\mXbar}_{1:2}}
      {1+n(n-1){\mXbar}'_{2}{\bf S}_{22}^{-1}{\mXbar}_{2}}}
\end{eqnarray}
and
\begin{eqnarray}
  U= n(n-1){\mXbar}'_{1:2}{\bf S}^{-1}_{11:2}{\mXbar}_{1:2}
\end{eqnarray}
respectively. In the literature, these two test statistics are called Rao's W and U statistics,
respectively. The test statistics $W$ is derived by the likelihood ratio principle.
Marden and Perlman [6] showed that for problem (1.2) both Rao's $W$-test and Rao's $U$-test are
similar and unbiased.

\indent  The invariance of problem (1.2) under a group ${\it G}$ of linear transformations,
where ${\it G}$ is the group of $p\times p$ nonsingular matrices of the form
\begin{align}
  {\bf g}=\left[
  \begin{array}{cc}
    {\bf g}_{11}& {\bf g}_{12}\\
    {\bf 0}& {\bf g}_{22}
  \end{array}
  \right]
\end{align}
with ${\bf g}_{11}:p_{1}\times p_{1}$ and ${\bf g}_{22}:p_{2}\times p_{2}$, can be exploited
so that the group ${\it G}$ acts on the sample space via ${\bf g}: ({\mXbar},~{\bf S})\to
({\bf g}{\mXbar}, ~{\bf g}{\bf S}{\bf g}^{'})$, and on the parameter space via
${\bf g}: ({\mbtheta},~{\mbSigma})\to ({\bf g}{\mbtheta},~{\bf g}{\mbSigma}{\bf g}^{'})$.
Let ${\mbtheta}_{1:2}$ and ${\mbSigma}_{11:2}$ be defined
similarly as in (1.5) and (1.6) but such that ${\mbtheta}$ and ${\mbSigma}$ replace ${\mXbar}$
and ${\bf S}$, respectively. Adopting the notions of Marden and Perlman [6], the maximal
invariant statistic is the pair $(L(1+M),~M)$ with $L=W$ defined as in (1.7) and
$M=n(n-1){\mXbar}'_{2}{\bf S}^{-1}_{22}{\mXbar}_{2}$, and correspondingly, the maximal
invariant parameter is the pair $({\Delta}_{1}, {\Delta}_{2})$ with
${\Delta}_{1}=n{\mbtheta}^{'}_{1:2}{\mbSigma}^{-1}_{11:2}{\mbtheta}_{1:2}$ and
${\Delta}_{2}=n{\mbtheta}^{'}_{2}{\mbSigma}^{-1}_{22}{\mbtheta}_{2}$.

\indent Using only ${\it G}$-invariant tests, Marden and Perlman ([6], p. 27)
concluded that the problem (1.2) reduces to that of testing
\begin{eqnarray}
 H^{I*}_{0}: {\Delta}= 0~\hbox{versus}~H^{I*}_{1}: {\Delta} > 0
\end{eqnarray}
based on $(L,~M)$, where ${\Delta}=n{\mbtheta}^{'}_{1}{\mbSigma}^{-1}_{11:2}{\mbtheta}_{1}$.
Marden and Perlman [6] established the necessary and sufficient conditions for the
admissibility of problem (1.10). They considered the homeomorphic transformations of
$(L, M)$, and inferred that Rao's $U$-test is inadmissible when only ${\it G}$-invariant
tests are applied to problem (1.2), but the overall Hotelling $T^{2}$-test is admissible.
However, their conclusions are against our statistical intuition. Also, their simulation
studies (see Tables 4.1a-4.1c) indicate a totally different story. As such, this work
attempts to clarify these phenomena.

\indent  First, note that if only the ${\it G}$-invariant tests are applied, then problem
(1.2) does not reduce to that of testing problem (1.10) but should reduce to that of testing
\begin{eqnarray}
H^{I}_{01}: {\Delta}_{1}= 0,~{\Delta}_{2}= 0 ~\hbox{versus}~H^{I}_{11}: {\Delta}_{1} >
0,~{\Delta}_{2}= 0.
\end{eqnarray}
Moreover, problems (1.11) and (1.10) are not identical. Problem (1.11) is easily seen
to imply problem (1.10), but not vice versa. Problem $(1.10)$ can be regarded as the union
of problem $(1.11)$ and the following subproblems: (i).
$H_{0}:{\Delta}_{1}= 0, ~{\Delta}_{2}= 0 ~
~\hbox{versus}~H_{1}:{\Delta}_{1} > 0,~{\Delta}_{2} > 0,$ (ii). $H_{0}:{\Delta}_{1}= 0,
~{\Delta}_{2} > 0 ~~\hbox{versus}~H_{1}:{\Delta}_{1} > 0,~{\Delta}_{2}= 0,$ and (iii).
$H_{0}: {\Delta}_{1}= 0,~{\Delta}_{2} > 0 ~~\hbox{versus}~H_{1}:{\Delta}_{1} > 0
,~{\Delta}_{2} > 0$. Problem $(1.11)$ provides more insight than does problem $(1.10)$ into
${\Delta}_{2}=0$ both in the null hypothesis and in the alternative hypothesis. In fact,
problem $(1.10)$ is a two-dimensional testing problem in which ${\Delta}_{2}$ is the nuisance
parameter, while problem $(1.11)$ is a one-dimensional testing problem. Problem $(1.11)$ is a
subproblem of problem $(1.10)$, so intuitively the optimal tests for problem $(1.11)$ (which is
equivalent to the hypothesis testing problem (1.2)) are not necessarily optimal for problem $(1.10)$,
and vice versa. Marden and Perlman [6] also clearly made this point (for details see the last three
lines of page 28 of Marden and Perlman [6]). Therefore, based on the optimal criterions set up for
the two-dimensional testing problem (1.10), to infer problem (1.11), which is only a one-dimensional
testing problem, the conclusions drawn may provide misleading messages. To ensure the information
$\mbmu_{2}={\bf 0}$ (i.e., ${\Delta}_{2}=0$) being incorporated for problem $(1.10)$, Marden and Perlman
[6] further made an assumption that $M$ is an ancillary statistic. Because the density function of
the ancillary statistic $M$ does not depend on the parameter ${\Delta}_{2}$ both under the null
hypothesis and under the alternative hypothesis, and hence the whole statistical inference should
depend on the conditional density function of $L$ given $M$, which is a noncentral $F$-type distribution,
in their set up. However, for the case of Hotelling's $T^{2}$-test, Marden and Perlman ([6], page 50)
adopted the exponential family for the problem (1.10) [not for the problem (1.11)] set up for their
statistical inference.

\indent Let ${\it Gl}$ be the general linear group. The problem of testing
$H^{u}_{0}:{\mbtheta} = {\bf 0}$ versus $H^{u}_{1}:{\mbtheta} \ne {\bf 0}$ is ${\it Gl}$-invariant.
When only ${\it Gl}$-invariant tests are performed, this problem reduces to that of testing
$H^{I}_{05}: {\Delta}^{*}=0,~\hbox{versus}~H^{I}_{15}: {\Delta}^{*} > 0$, where
${\Delta}^{*}={\Delta}_{1}+{\Delta}_{2}$. For this ${\it Gl}$-invariant testing problem,
Hotelling's $T^{2}$-test is well-known to be the uniformly most powerful test
(Simaika [14]), and so is admissible. Schwartz [13] employed the Birnbaum-Stein method (Birnbaum
[2], Stein [15]) to study the admissibility of fully ${\it Gl}$-invariant tests in the
multivariate analysis of variance setting. Problem (1.2) is not ${\it Gl}$-invariant,
although it is ${\it G}$-invariant. Therefore, studying the power domination problems
of Hotelling's $T^{2}$-test, Rao's $W$-test and Rao's $U$-test for problem (1.2) via the fully
${\it Gl}$-invariant approach may not be helpful. The group ${\it G}$ is amenable and
meets the conditions of the Hunt-Stein theorem (Lehmann [5]). Therefore, any minimax questions
in problem (1.2) can be reduced by the group ${\it G}$.



\indent Notably, problems (1.2) and (1.11) are the problems of testing against restricted
alternatives. However, problem (1.10) is not such a problem. Therefore, by utilizing problem (1.10)
to draw inferences concerning problem (1.2), we may overlook the intrinsic nature of the restricted
alternative (because $\Delta_2=0$, which is determined directly from the basic assumption
$\mbtheta_2={\bf 0})$ when applying results in the literature or developing new theories.

\indent 
The exponential
structure of the distribution of $(\mXbar, {\bf S})$ is incorporated to generalize the Birnbaum-Stein
method for problem (1.2). As a result, Section 2 presents two main results: the acceptance region of
Rao's $U$-test is convex, and Rao's $U$-test is admissible. Section 3 applies Eaton's [3] results
to show that Hotelling's $T^{2}$-test is inadmissible for problem (1.2). The discussion regarding
the Rao $W$-test is given in Section 4. Some general remarks are given in the final section.
The Appendix provides six lemmas to show the convexity of the acceptance region of Rao's $U$-test.

\vspace{0.3cm}
\def \theequation{2.\arabic{equation}}
\setcounter{equation}{0}
\noindent {\bf 2. Admissibility of Rao's $U$-test}
\vspace{0.3cm}

\indent Stein [15] proved in detail that Hotelling's $T^{2}$-test is admissible for
the problem of testing $H^{u}_{0}:{\mbtheta} = {\bf 0}$ against the global alternative
$H^{u}_{1}:{\mbtheta} \ne {\bf 0}$. This proof can also be found in Anderson's book
([1], p. 188-190). The main purpose of this section is to incorporate the Birnbaum-Stein method
to demonstrate that Rao's $U$-test is admissible for problem (1.2). Whether the acceptance region of
Rao's $U$-test is a convex set must first be determined. Let
\begin{eqnarray}
{\mathcal A}_{U} = \{\, (\mXbar, \mbS) \mid n(n-1)\mXbar_{1:2}'\mbS_{11:2}^{-1}
\mXbar_{1:2} \leq k,~{\mbS} \mbox{ is p.d.} \,\}
\end{eqnarray}
for a suitable k, and

\begin{eqnarray}
\pmb{B}^{+}(\mbS)
&=& \left(
\begin{array}{c}
{\bf I}\\
-\mbS_{22}^{-1}\mbS_{21}
\end{array}
\right)
\mbS_{11:2}^{-1}
\left(
\begin{array}{cc}
{\bf I}& -\mbS_{12}\mbS_{22}^{-1}
\end{array}
\right) \\ \nonumber
&=& \mbS^{-1}
- \left[
\begin{array}{cc}
{\bf 0}& {\bf 0}\\
{\bf 0}& \mbS_{22}^{-1}
\end{array}
\right]_,
\end{eqnarray}
which is positive semi-definite (p.s.d.) and of rank $p_1$. Also let

\begin{align}
{\bf B}(\mbS)= \left(
\begin{array}{c}
{\bf I}\\
-\mbS_{22}^{-1}\mbS_{21}
\end{array}
\right)
({\bf I} + \mbS_{12}\mbS_{22}^{-1}\mbS_{22}^{-1}\mbS_{21})^{-1}\mbS_{11:2}
({\bf I} + \mbS_{12}\mbS_{22}^{-1}\mbS_{22}^{-1}\mbS_{21})^{-1} \left(
\begin{array}{cc}
{\bf I}& -\mbS_{12}\mbS_{22}^{-1}
\end{array}
\right)_.
\end{align}
Let ${\bf C}^{+}$ be the Moore-Penrose generalized inverse of ${\bf C}$.
Then,
\begin{eqnarray}
{\bf B}(\mbS) = \left(
\begin{array}{cc}
{\bf I}& -\mbS_{12}\mbS_{22}^{-1}
\end{array}
\right)^{+} {\mbS}_{11:2} \left(
\begin{array}{c}
{\bf I}\\
-\mbS_{22}^{-1}\mbS_{21}
\end{array}
\right)^{+}_,
\end{eqnarray}
which is shown to be the Moore-Penrose generalized inverse of ${\bf B}^{+}(\mbS)$ in
Lemma 2 of the Appendix, is easily established. For the notions related to the Moore-Penrose
generalized inverse of matrices and matrix-convex (matrix-concave) functions, refer to
Rao and Mitra [11] and Marshall and Olkin [7], respectively. The Appendix develops
six lemmas related to the Moore-Penrose generalized inverse of matrices, matrix-convex
and matrix-concave.

\indent Anderson ([1], problem 17 of page 193) claimed that the acceptance region of
Hotelling $T^{2}$-test
\begin{eqnarray}
{\mathcal A}_{T^{2}} = \{\, (\mXbar, \mbS) \mid n(n-1)\mXbar'{\mbS}^{-1}\mXbar \leq k_{1} \}
~~\hbox{is a convex set}.
\end{eqnarray}
Notably, by (2.1) and (2.2) the accepted region of Rao's $U$-test ${\mathcal A}_{U}$ can be rewritten as
\begin{eqnarray}
{\mathcal A}_{U} = \{\, (\mXbar, \mbS) \mid n(n-1)\mXbar'{\bf B}^{+}(\mbS)
\mXbar \leq k \},
\end{eqnarray}
where ${\bf B}^{+}(\mbS)$ is defined in (2.2). Let ${\bf A}\succeq {\bf B}$ denote that
the matrix ${\bf A}-{\bf B}$ is p.s.d. throughout this paper. The lemmas developed
in the Appendix are used to generalize Anderson's result (2.5) to the following theorem.

\vspace{0.3cm}
\noindent {\bf Theorem 1.} {\it Let ${\cal S}$ be the set of all $p\times p$ positive definite
matrices. Then ${\mathcal A}_{U} = \{\, (\mXbar, \mbS) \mid n(n-1)\mXbar_{1:2}'\mbS_{11:2}^{-1}
\mXbar_{1:2} \leq k \,\}$ is convex on $R^{p}\times {\cal S}$}.

\indent {\bf Proof.} By (2.6), ${\mathcal A}_{U} = \{\, (\mXbar, \mbS) \mid n(n-1)\mXbar'{\bf B}^{+}
(\mbS)\mXbar \leq k \,\}$. Let ${\bf B}(\mbS)$ be defined as in (2.3), then by Lemma 2
in the Appendix, it is the Moore-Penrose generalized inverse of ${\bf B}^{+}(\mbS)$. Lemma
$4$ shows that ${\bf B}(\mbS)$ is matrix concave on ${\cal S}$,
that is, $\forall\, \alpha \in (0, 1)$
\begin{align}
{\bf B}(\alpha\mbS + (1 - \alpha)\mbT) \succeq \alpha{\bf B}(\mbS)
+ (1 - \alpha){\bf B}(\mbT).
\end{align}
Therefore, by Lemma 5,
\begin{eqnarray}
{\bf B}^{+}(\alpha\mbS + (1 - \alpha)\mbT) \preceq
\left( \alpha{\bf B}(\mbS) + (1 - \alpha){\bf B}(\mbT) \right)^{+}.
\end{eqnarray}
By the inequality (2.8) and Lemma 6,
\begin{align}
&(\alpha\mXbar + (1 - \alpha)\mYbar)'{\bf B}^{+}(\alpha\mbS + (1 - \alpha)\mbT)
(\alpha\mXbar + (1 - \alpha)\mYbar)  \\ \nonumber
&\qquad \leq (\alpha\mXbar + (1 - \alpha)\mYbar)'(\alpha{\bf B}(\mbS) + (1 - \alpha)
{\bf B}(\mbT))^{+}(\alpha\mXbar + (1 - \alpha)\mYbar)  \\ \nonumber
&\qquad \leq \alpha\mXbar'{\bf B}^{+}(\mbS)\mXbar + (1 - \alpha)\mYbar'{\bf B}^{+}
(\mbT)\mYbar.
\end{align}

\indent For the definitions of $\alpha$-admissible and $d$-admissible, we may refer to the page 306 of
Lehmann [6].

\indent {\bf Remark}. For testing against global alternative, (i.e., $H_{0}: {\mbtheta}={\bf 0}$
against $H_{1}: {\mbtheta}\ne {\bf 0}$), since the Hotelling's $T^{2}$-test statistic
$n(n-1)\Bar{\bf X}^{'}{\bf S}^{-1}\Bar{\bf X}$ is the monotone function of
$n(n-1)\Bar{\bf X}^{'}{\bf W}^{-1}\Bar{\bf X}$, where ${\bf W}={\bf S}+n\Bar{\bf X}\Bar{\bf X}^{'}$,
Stein worked on the space $(\Bar{\bf X}, {\bf W})$. Theorem 5.6.6 of Anderson [1] does not require that it
is necessary to work on the space $(\Bar{\bf X}, {\bf W})$. Note that for the testing
$H_{0}: {\mbtheta}_{1}={\bf 0}, {\mbtheta}_{2}={\bf 0}$ against
the restricted alternative $H_{1}: {\mbtheta}_{1}\ne {\bf 0}, {\mbtheta}_{2}={\bf 0}$, the Rao's $U$-test
statistic $n(n-1)\Bar{\bf X}^{'}\pmb{B}^{+}(\mbS)\Bar{\bf X}$ is no longer to be a monotone function of
$n(n-1)\Bar{\bf X}^{'}\pmb{B}^{+}({\bf W})\Bar{\bf X}$ any more. To overcome the difficulty, we work
on the space $(\Bar{\bf X}, {\bf S})$ instead. Let ${\mathcal H}^{*}_{c} =\{\, (\mXbar, \mbS)
\mid n{\mbtheta}^{'}{\mbSigma}^{-1}\mXbar -\frac{1}{2} \tr({\mbSigma}^{-1}\mbS)  > c \,\}$.  It is easy to note
that ${\mathcal H}^{*}_{c} $ is a half-space, and hence the assumption that the acceptance region $A_{U}$ is
disjoint with the half-space ${\mathcal H}^{*}_{c}$ on the space $(\Bar{\bf X}, {\bf S})$ holds. Let
\begin{eqnarray}
{\mathcal H}^{a}_{c} = \{\, (\mXbar, \mbS) \mid n\mbtheta^{'}{\mbSigma}^{-1}\mXbar
-\frac{1}{2} \tr [{\mbSigma}^{-1}(\mbS + n\mXbar\mXbar')] > c \,\}.
\end{eqnarray}
Note that
$(\mXbar, \mbS) \in {\mathcal H}^{a}_{c}$ implies that $(\mXbar, \mbS) \in {\mathcal H}^{*}_{c}$. Thus,
\begin{eqnarray}
{\mathcal H}^{a}_{c} \subset {\mathcal H}^{*}_{c}.
\end{eqnarray}
And hence, the intersection of ${\mathcal A}_{U}$ and ${\mathcal H}^{a}_{c}$ is also empty. As such, we may
work the proof of Theorem 2 on the space $(\Bar{\bf X}, {\bf S})$.

\vspace{0.3cm}
\noindent {\bf Theorem 2.} {\it For the problem (1.2), Rao's $U$-test is $\alpha$-admissible}.

\indent {\bf Proof.} The likelihood function of $\boldsymbol{X}_1$, $\cdots$, $\boldsymbol{X}_n$ is
\begin{align}
\frac{e^{-\frac{1}{2}n\mbtheta'\mbSigma^{-1}\mbtheta}}{(2\pi)^{\frac{1}{2}pn}|\mbSigma|
^{\frac{n}{2}}}\exp[n\mbtheta'\mbSigma^{-1}\mXbar + \tr(-\frac{1}{2}
\mbSigma^{-1}\sum_{i=1}^n{\bf X}_i{\bf X}_i')].
\end{align}
Let $\boldsymbol{w} = (\mbw^{(1)'}, \mbw^{(2)'})'$, where
$\mbw^{(1)} = \mbSigma^{-1}\mbtheta$ and
$\mbw^{(2)} = -\frac{1}{2}(\sigma^{11}, \cdots, \sigma^{1p},
\sigma^{22}, \cdots, \sigma^{pp})'$,
where $(\sigma^{ij}) = \mbSigma^{-1}$. Let $\mbw^{(1)'} =\mbnu' = (\mbnu_1', \mbnu_2')'$.
By Theorem 1, ${\mathcal A}_{U}$ is convex on $R^{p}\times {\cal S}$. Consider the other condition
of theorem 5.6.5. (Anderson [1]), ${\mathcal A}_{U}$ is assumed to be disjoint with the subspace
\begin{eqnarray}
{\mathcal H}_{c} = \{\, (\mXbar, \mbS) \mid n\mbnu'\mXbar
-\frac{1}{2} \tr\mbLambda(\mbS + n\mXbar\mXbar') > c \,\},
\end{eqnarray}
where $\mbLambda$ is symmetric, for some $c$.

\indent Theorem 8 presented by Lehmann ([5], page 307) can be applied if $\mbw_0 + \lambda\mbw_1 \in
H_1$ can be demonstrated for $\lambda > 0$, which can be accomplished with the following two steps: (I)
${\bf I} + \lambda\mbLambda$ is p.d., and
(I\!I) $\mbtheta_0 + \lambda\mbtheta \in H_1$ for $\lambda > 0$.

(I). That ${\bf I} + \lambda\mbLambda$ is p.d. is shown if $\mbLambda$ can be shown to be
p.s.d., for $\lambda > 0$. Suppose that $\mbLambda$ is not p.s.d., then by arguments similar
to those in Anderson ([1], p. 189-190) it can be written
\[
{\mbLambda}= {\bf D}\left[\begin{array}{cccc}
    {\bf I}   &  {\bf 0}  &  {\bf 0}  \\
    {\bf 0}   & -{\bf I}  &  {\bf 0}  \\
    {\bf 0}   &  {\bf 0}  &  {\bf 0}  \\
\end{array} \right]{\bf D}^{'},\]  \\
where ${\bf D}$ is nonsingular. Let $\mXbar=(1/{\gamma}){\bf X}_{0}$ and
\[
{\bf S}= ({\bf D}^{'})^{-1}\left[\begin{array}{cccc}
    {\bf I}   &  {\bf 0}  &  {\bf 0}  \\
    {\bf 0}   &  {\gamma}{\bf I}  &  {\bf 0}  \\
    {\bf 0}   &  {\bf 0}  &  {\bf I}  \\
\end{array} \right]{\bf D}^{-1}.
\]
Then
\begin{eqnarray}
&& n\mbnu'\mXbar-\frac{1}{2} \tr\mbLambda(\mbS + n\mXbar\mXbar')  \\ \nonumber
&=& \frac{n}{\gamma}{\mbnu}^{'}{\bf X}_{0}-\frac{n}{2{\gamma}^{2}}{\bf X}^{'}_{0}
{\bf D}\left[\begin{array}{cccc}
    {\bf I}   &  {\bf 0}  &  {\bf 0}  \\
    {\bf 0}   & -{\bf I}  &  {\bf 0}  \\
    {\bf 0}   &  {\bf 0}  &  {\bf 0}  \\
\end{array} \right]{\bf D}^{'}{\bf X}_{0}
+\frac{1}{2}\tr\left[\begin{array}{cccc}
    -{\bf I}   &  {\bf 0}  &  {\bf 0}  \\
    {\bf 0}   &  {\gamma}{\bf I}  &  {\bf 0}  \\
    {\bf 0}   &  {\bf 0}  &  {\bf I}  \\
\end{array} \right],  \\  \nonumber
\end{eqnarray}
which is greater than c for sufficiently large $\gamma$. Thus, the subspace ${\mathcal H}_{c}$
reduces to
\begin{eqnarray}
{\mathcal H}^{\gamma}_{c} &=& \{\, {\bf X}_{0}  \mid
\frac{n}{\gamma}{\mbnu}^{'}{\bf X}_{0}-\frac{n}{2{\gamma}^{2}}{\bf X}^{'}_{0}
{\bf D}\left[\begin{array}{cccc}
    {\bf I}   &  {\bf 0}  &  {\bf 0}  \\
    {\bf 0}   & -{\bf I}  &  {\bf 0}  \\
    {\bf 0}   &  {\bf 0}  &  {\bf 0}  \\
\end{array} \right]{\bf D}^{'}{\bf X}_{0}   \\  \nonumber
& & \hspace{5cm} +\frac{1}{2}\tr\left[\begin{array}{cccc}
    -{\bf I}   &  {\bf 0}  &  {\bf 0}  \\
    {\bf 0}   &  {\gamma}{\bf I}  &  {\bf 0}  \\
    {\bf 0}   &  {\bf 0}  &  {\bf I}  \\
\end{array} \right] > c \,\}
\end{eqnarray}
for sufficiently large $\gamma$. Obviously, ${\mathcal H}^{\gamma}_{c}=R^{p}$ as $\gamma$
approaches infinity.  Now, let
\begin{align}
{\mathcal A}_{T^2, k^{\star}} &= \{\, (\mXbar, \mbS) \mid n(n-1)\mXbar'
    (\mbS + n\mXbar\mXbar')^{-1}\mXbar \leq k^{\star} \,\} \\ \nonumber
              &= \{\, (\mXbar, \mbS) \mid n(n-1)\mXbar'\mbS^{-1}\mXbar \leq k \,\}
\end{align}
with $k=k^{*}/(1-k^{*})$. Then, (2.16) reduces to
\begin{eqnarray}
  {\mathcal A}^{\gamma}_{T^2, k^{*}}=\{\, {\bf X}_{0} \mid
   \frac{n(n-1)}{{\gamma}^{2}}{\bf X}^{'}_{0}
{\bf D}\left[\begin{array}{cccc}
    {\bf I}   &  {\bf 0}  &  {\bf 0}  \\
    {\bf 0}   &  {\gamma}^{-1}{\bf I}  &  {\bf 0}  \\
    {\bf 0}   &  {\bf 0}  &  {\bf I}  \\
\end{array} \right]{\bf D}^{'}{\bf X}_{0} \leq k\},  \\  \nonumber
\end{eqnarray}
for sufficiently large $\gamma$. It can be easily seen that ${\mathcal H}^{\gamma}_{c}\cap
{\mathcal A}^{\gamma}_{T^2, k^{*}} \ne \varnothing$  for sufficiently large $\gamma$.
Furthermore,
\begin{align}
{\mathcal A}_{T^2, k^{*}} &= \{\, (\mXbar, \mbS) \mid n(n-1)\mXbar'\mbS^{-1}\mXbar
            \leq k \,\} \\  \nonumber
            &= \{\, (\mXbar, \mbS) \mid n(n-1)\mXbar_{1:2}'\mbS_{11:2}^{-1}\mXbar_{1:2}
             + n(n-1)\mXbar_{2}'\mbS_{22}^{-1}\mXbar_2 \leq k \,\}  \\ \nonumber
            &\subseteqq \{\, (\mXbar, \mbS) \mid n(n-1)\mXbar_{1:2}'\mbS_{11:2}^{-1}
            \mXbar_{1:2} \leq k \,\}  \\  \nonumber
            &= {\mathcal A}_{U}.
\end{align}
Thus, if $\pmb{\Lambda}$ is not p.s.d., then
\begin{align}
{\mathcal A}_{U} \cap \mathcal{H}_{c} \neq \varnothing,
\end{align}
which leads to a contradiction. Therefore, $\mbLambda$ is p.s.d..

\indent To proceed with step (I\!I),
note that ${\bf I}+{\mbLambda}$ is p.d., and without loss of generality, its
inverse can be denoted by ${\mbSigma}$, so $({\bf I}+{\mbLambda})^{-1}={\mbSigma}$. Some
more notation is needed. Let
\begin{align}
\mbg = \left[
\begin{array}{cc}
{\bf I}& -\mbSigma_{12}\mbSigma_{22}^{-1}\\
{\bf 0}& {\bf I}
\end{array}
\right], \quad
(\mbtheta, \mbSigma) \stackrel{\mbg}{\longrightarrow}
(\mbg\mbtheta, \mbg\mbSigma\mbg')
\end{align}
and write
\begin{align}
\widetilde{\mbSigma} &= \mbg \mbSigma \mbg'   \\  \nonumber
&= \left[
\begin{array}{cc}
\mbSigma_{11} - \mbSigma_{12} \mbSigma_{22}^{-1} \mbSigma_{21}& {\bf 0}\\
{\bf 0}& \mbSigma_{22}
\end{array}
\right]  \\  \nonumber
&= \left[
\begin{array}{cc}
\mbSigma_{11:2}& {\bf 0}\\
{\bf 0}& \mbSigma_{22}
\end{array}
\right]_.
\end{align}
Notably, $\widetilde{\mbSigma} = \widetilde{\mbSigma}'$.
Let $\mbZ = \mbg \mXbar$ and $\mbS_0=\mbg \mbS \mbg'$.
Then,
\begin{align}
\mbZ = \left(
\begin{array}{c}
\mXbar_1 - \mbSigma_{12}\mbSigma_{22}^{-1}\mXbar_2\\
\mXbar_2
\end{array}
\right)_{,}
\end{align}
\begin{align}
\mbS_0
= \left[
\begin{array}{cc}
\mbS_{11} - \mbSigma_{12}\mbSigma_{22}^{-1}\mbS_{21}
- \mbS_{12}\mbSigma_{22}^{-1}\mbSigma_{21}
+ \mbSigma_{12}\mbSigma_{22}^{-1}\mbS_{22}\mbSigma_{22}^{-1}\mbSigma_{21}&
\mbS_{12} - \mbSigma_{12}\mbSigma_{22}^{-1}\mbS_{22}\\
\mbS_{21} - \mbS_{22}\mbSigma_{22}^{-1}\mbSigma_{21}& \mbS_{22}
\end{array}
\right]_,
\end{align}
and
\begin{align}
&\tr \mbSigma^{-1}(\mbS + n\mXbar\mXbar')   \\  \nonumber
&= \tr (\mbg \mbSigma \mbg')^{-1}\mbg
(\mbS + n\mXbar\mXbar')\mbg'  \\  \nonumber
&= \tr \widetilde{\mbSigma}^{-1}(\mbS_0 + n\mbZ\mbZ')_,
\end{align}
\begin{align}
 {\mbnu}^{'}\mXbar &=\widetilde{\mbnu}'\mbZ, \mbox{ where } \\ \nonumber
 \hspace{1cm} \widetilde{\mbnu}
&  =({\mbg}^{-1})^{'} \mbnu      \\  \nonumber
& = \left(
\begin{array}{c}
{\mbSigma}^{-1}_{11:2}({\mbtheta}_{1}-{\mbSigma}_{12}{\mbSigma}^{-1}_{22}{\mbtheta}_{2})\\
{\mbSigma}^{-1}_{22}{\mbtheta}_{2}
\end{array}
\right)   \\  \nonumber
&\triangleq \left(
\begin{array}{c}
\widetilde{\mbnu}_{1}\\
\widetilde{\mbnu}_{2}
\end{array}
\right).
\end{align}
Thus, the subspace ${\mathcal H}^{a}_{c}$ becomes
\begin{align}
{\mathcal H}^{a}_{c} = \{\, (\mXbar, \mbS) \mid n\widetilde{\mbnu}'\mbZ
- \frac{1}{2} \tr \widetilde{\mbSigma}^{-1}(\mbS_0 + n\mbZ\mbZ') > c \}_.
\end{align}

\indent Equation (2.25) indicates that the hypothesis testing problem (1.2)
$H_{0}: {\mbtheta}_{1}={\bf 0}, {\mbtheta}_{2}={\bf 0}$ versus
$H_{1}: {\mbtheta}_{1}\ne{\bf 0}, {\mbtheta}_{2}={\bf 0}$ is equivalent to the hypothesis
testing problem
\begin{eqnarray}
H^{*}_{0}: \widetilde{\mbnu}_{1}={\bf 0}, \widetilde{\mbnu}_{2}={\bf 0}~~\mbox{versus}~~
H^{*}_{1}: \widetilde{\mbnu}_{1}\ne{\bf 0}, \widetilde{\mbnu}_{2}={\bf 0}.
\end{eqnarray}
And hence to show that $\mbtheta_0 + \lambda\mbtheta \in H_1$, $\lambda > 0$ for problem (1.2)
is equivalent to showing that
$\widetilde\mbnu_0 + \lambda\widetilde\mbnu \in H_1$, $\lambda > 0$ for problem (2.27).
Next, step (I\!I) is considered.

\indent (I\!I). To show that $\widetilde\mbnu_0 + \lambda\widetilde\mbnu \in H_1$ for
$\lambda > 0$, the aim is to demonstrate that $\widetilde{\mbnu}_1 \neq {\bf 0}$ and
$\widetilde{\mbnu}_2 = {\bf 0}$. If the statement that
$\widetilde{\mbnu}_1 \neq {\bf 0}$ and $\widetilde{\mbnu}_2 = {\bf 0}$ is not true,
then, there are three cases (i) $\widetilde{\mbnu}_1 \neq {\bf 0}$,
$\widetilde{\mbnu}_2 \neq {\bf 0}$, (ii) $\widetilde{\mbnu}_1 = {\bf 0}$,
$\widetilde{\mbnu}_2 \neq {\bf 0}$ and (iii) $\widetilde{\mbnu}_1 = {\bf 0}$,
$\widetilde{\mbnu}_2 = {\bf 0}$. We assume that (i), (ii) and (iii) are true, and then show that
those assumptions to lead to contradictions. To proceed, it is enough to consider the situation that
${\bf S}={\mbSigma}$.

\indent Note that given ${\bf S}={\mbSigma}$, then
\begin{eqnarray}
\mbS_0
= \left[
\begin{array}{cc}
\mbSigma_{11:2}& {\bf 0}\\
{\bf 0}& \mbSigma_{22}
\end{array}
\right]_.
\end{eqnarray}
Accordingly, both sets ${\mathcal A}_{U}$ and ${\mathcal H}^{a}_{c}$ are reduced to
p-dimensional sets (using the notation loosely)
\begin{align}
{\mathcal A}_{U} = \{\, \mbZ \mid n(n-1) \mbZ_1'\mbSigma_{11:2}^{-1}\mbZ_1 \leq k \,\}
\end{align}
and
\begin{align}
{\mathcal H}^{a}_c = \{\, \mbZ \mid n\widetilde{\mbnu}'\mbZ - \frac{p}{2} -
\frac{n}{2}\mbZ_1'\mbSigma_{11:2}^{-1}\mbZ_1 - \frac{n}{2}\mbZ_2'\mbSigma_{22}^{-1}
\mbZ_2 > c \,\},
\end{align}
respectively. Notably, given that $\mbS = \mbSigma$, the assumption that the sets
${\mathcal A}_{U}$ and ${\mathcal H}^{a}_c$ are disjoint still holds.

\indent First, (i) $\widetilde{\mbnu}_1 \neq {\bf 0}$, $\widetilde{\mbnu}_2 \neq {\bf 0}$
is assumed. In the problem of testing $H_0^1$: $\widetilde{\mbnu}_1 = {\bf 0}$ versus
$H_1^1$: $\widetilde{\mbnu}_1 \neq {\bf 0}$, whenever $\widetilde{\mbnu}_1 \neq {\bf 0}$,
then there exists a constant $c_{1}$ which does not depend on $\mbZ_1$ and $\mbSigma_{11:2}$
such that
\begin{align}
{\mathcal H}_{c_1}^{a\star} \cap {\mathcal A}_{T^2, k_1^{\star}}
= \varnothing \quad  \mbox{and} \quad {\mathcal H}_{c_1 - \epsilon}^{a\star}
 \cap {\mathcal A}_{T^2,k_1^{\star}} \neq \varnothing
\end{align}
for any $\epsilon > 0$, where
\begin{align}
{\mathcal H}_{c_1}^{a\star} = \{\, \mbZ_1 \mid n\widetilde{\mbnu}_1'\mbZ_1
 - \frac{p_1}{2} - \frac{n}{2}\mbZ_1'\mbSigma_{11:2}^{-1}\mbZ_1 > c_1+\frac{p_2}{2} \,\}
\end{align}
and
\begin{align}
{\mathcal A}_{T^2, k_1^{\star}}
&= \{\, \mbZ_1 \mid n (n-1) \mbZ_1'(\mbSigma_{11:2} + n\mbZ_1\mbZ_1')^{-1}\mbZ_1
    \leq k_1^*\} \\ \nonumber
&= \{\, \mbZ_1 \mid n(n-1) \mbZ_1'\mbSigma_{11:2}^{-1}\mbZ_1 \leq k_1 \,\}.
\end{align}
This is equivalent to the existence of a constant $c$ such that
\begin{eqnarray}
{\mathcal H}_c^{a\star} \cap {\mathcal A}_{T^2, k^{\star}}
= \varnothing \quad \mbox{and} \quad {\mathcal H}_{c - \epsilon}^{a\star}
\cap {\mathcal A}_{T^2, k^{\star}} \neq \varnothing
\end{eqnarray}
for any $\epsilon > 0$. Rewrite ${\mathcal H}^{a}_c$ as
\begin{align}
{\mathcal H}^{a}_c = \{\, \mbZ \mid n\widetilde{\mbnu}_1'\mbZ_1 - \frac{p}{2}
- \frac{n}{2}\mbZ_1'\mbSigma_{11:2}^{-1}\mbZ_1
+ \frac{n}{2}[\widetilde{\mbnu}_2'\mbSigma_{22}\widetilde{\mbnu}_{2}
- (\mbZ_2 - \mbSigma_{22}\widetilde{\mbnu}_2)'\mbSigma_{22}^{-1}
(\mbZ_2 - \mbSigma_{22}\widetilde{\mbnu}_2)] > c \,\}.
\end{align}
Consider
\begin{eqnarray}
\epsilon ({\bf Z}_{2}) =
\frac{n}{2}[\widetilde{\mbnu}'_{2}\mbSigma_{22}\widetilde{\mbnu}_{2}
- (\mbZ_2 -\mbSigma_{22}\widetilde{\mbnu}_2)'\mbSigma_{22}^{-1}
(\mbZ_2 - \mbSigma_{22}\widetilde{\mbnu}_2)]~>~0.
\end{eqnarray}
Then, from equation (2.34),
\begin{eqnarray}
{\mathcal H}^{a\star}_{c - \epsilon ({\bf Z}_{2})}
\cap {\mathcal A}_{T^2, k^{\star}} \neq \varnothing,
\end{eqnarray}
where
\begin{eqnarray}
{\mathcal H}^{a\star}_{c - \epsilon ({\bf Z}_{2})}
&=& \{\, {\mbZ}_{1} \mid n\widetilde{\mbnu}'_{1}\mbZ_1 -
\frac{p}{2} - \frac{n}{2}\mbZ_1'\mbSigma_{11:2}^{-1}\mbZ_1   \\ \nonumber
& &+ \frac{n}{2}[\widetilde{\mbnu}'_{2}\mbSigma_{22}\widetilde{\mbnu}_{2}
- (\mbZ_2 - \mbSigma_{22}\widetilde{\mbnu}_2)'\mbSigma_{22}^{-1}
(\mbZ_2 - \mbSigma_{22}\widetilde{\mbnu}_2)] > c \,\}.
\end{eqnarray}
Hence,
\begin{eqnarray}
{\mathcal H}^{a\star}_{c - \epsilon ({\bf Z}_{2})}
\cap ({\mathcal A}_{T^2, k^{\star}} \times \mathbb{R}^{p_2}) \neq \varnothing.
\end{eqnarray}
Let
\begin{eqnarray}
{\mathcal B}(\mbSigma_{22}\widetilde{\mbnu}_2) = \{\, \mbZ_2 \mid
(\mbZ_2 - \mbSigma_{22}\widetilde{\mbnu}_2)'\mbSigma_{22}^{-1}
(\mbZ_2 - \mbSigma_{22}\widetilde{\mbnu}_2)
\leq \widetilde{\mbnu}_2'\mbSigma_{22}\widetilde{\mbnu}_2 \,\}.
\end{eqnarray}
$\mbSigma_{22}\widetilde{\mbnu}_2 \in {\mathcal B}(\mbSigma_{22}\widetilde{\mbnu}_2)$,
so ${\mathcal B}(\mbSigma_{22}\widetilde{\mbnu}_2) \neq \varnothing$; and
$\mbZ_2 \in {\mathcal B}(\mbSigma_{22}\widetilde{\mbnu}_2)$ implies that
$\epsilon ({\bf Z}_{2})~>~0$.
Let
\begin{eqnarray}
 {\mathcal S}_{c}=\cup_{{\bf Z}_{2} \in {\mathcal B}(\mbSigma_{22}\widetilde{\mbnu}_2) }
  {\mathcal H}^{a\star}_{c - \epsilon ({\bf Z}_{2})}.
\end{eqnarray}
Now, (2.39) implies
\begin{eqnarray}
  {\mathcal S}_{c} \cap ({\mathcal A}_{T^2, k^{\star}}
   \times \mathbb{R}^{p_2}) \neq \varnothing.
\end{eqnarray}
Notably,
\begin{eqnarray}
 {\mathcal S}_{c}&=&\cup_{{\bf Z}_{2} \in {\mathcal B}(\mbSigma_{22}\widetilde{\mbnu}_2) }
  \{\, {\mbZ}_{1} \mid n\widetilde{\mbnu}_1'\mbZ_1 -
\frac{p}{2} - \frac{n}{2}\mbZ_1'\mbSigma_{11:2}^{-1}\mbZ_1      \\   \nonumber
& &~~~~~~~~~~~~~~~~~~~~~~~\vspace{2cm}+
\frac{n}{2}[\widetilde{\mbnu}'_{2}\mbSigma_{22}\widetilde{\mbnu}_{2} - (\mbZ_2 -
\mbSigma_{22}\widetilde{\mbnu}_2)'\mbSigma_{22}^{-1}
(\mbZ_2 - \mbSigma_{22}\widetilde{\mbnu}_2)] >
c \,\} \\   \nonumber    &=&\{\, {\mbZ} \mid n\widetilde{\mbnu}_1'\mbZ_1 -
\frac{p}{2} - \frac{n}{2}\mbZ_1'\mbSigma_{11:2}^{-1}\mbZ_1
+ \frac{n}{2}[\widetilde{\mbnu}_2'\mbSigma_{22}\widetilde{\mbnu}_{2}
- (\mbZ_2 - \mbSigma_{22}\widetilde{\mbnu}_2)'\mbSigma_{22}^{-1}
(\mbZ_2 - \mbSigma_{22}\widetilde{\mbnu}_2)] > c \,  \\   \nonumber
 & & \quad \vspace{2cm}~ \quad  ~~~~~~~~~~~~~~~~~~~~~~~~~\mbox{and}~ \quad
  \widetilde{\mbnu}_2'\mbGamma_{22}\widetilde{\mbnu}_{2}
\geq (\mbZ_2 - \mbSigma_{22}\widetilde{\mbnu}_2)'\mbSigma_{22}^{-1}
(\mbZ_2 - \mbSigma_{22}\widetilde{\mbnu}_2)\}  \\ \nonumber
& \subseteq & \{\, {\mbZ} \mid n\widetilde{\mbnu}_1'\mbZ_1 -
\frac{p}{2} - \frac{n}{2}\mbZ_1'\mbSigma_{11:2}^{-1}\mbZ_1
+ \frac{n}{2}[\widetilde{\mbnu}_2'\mbSigma_{22}\widetilde{\mbnu}_{2}
- (\mbZ_2 - \mbSigma_{22}\widetilde{\mbnu}_2)'\mbSigma_{22}^{-1}
(\mbZ_2 - \mbSigma_{22}\widetilde{\mbnu}_2)] > c \,\}  \\ \nonumber
&=& {\mathcal H}^{a}_{c}.
\end{eqnarray}
Thus,
\begin{eqnarray}
{\mathcal H}^{a}_c \cap ({\mathcal A}_{T^2, k^{\star}}
  \times \mathbb{R}^{p_2}) \neq \varnothing.
\end{eqnarray}
Notably,
\begin{align}
{\mathcal A}_{T^2, k^{\star}} \times \mathbb{R}^{p^2}
= \{\, \mbZ \mid n(n-1) \mbZ_1'\mbSigma_{11:2}^{-1}\mbZ_1 \leq k \,\}
={\mathcal A}_{U}.
\end{align}
Namely, if $\widetilde{\mbnu}_1 \neq {\bf 0}$ and $\widetilde{\mbnu}_2 \neq {\bf 0}$,
then ${\mathcal H}^{a}_c \cap {\mathcal A}_{U} \neq \varnothing$. This implies that
${\mathcal H}_c \cap {\mathcal A}_{U} \neq \varnothing$, and leads to a contradiction.

\indent Next, (ii) $\widetilde{\mbnu}_1 = {\bf 0}$, $\widetilde{\mbnu}_2 \neq {\bf 0}$
is assumed. Then, the set ${\mathcal H}^{a}_c$ in equation (2.35) reduces to
\begin{eqnarray}
{\mathcal H}^{a}_c = \{\, \mbZ
\mid n\mbZ_1'\mbSigma_{11:2}^{-1}\mbZ_1
+ n(\mbZ_2 - \mbSigma_{22}\widetilde{\mbnu}_2)'\mbSigma_{22}^{-1}
(\mbZ_2 - \mbSigma_{22}\widetilde{\mbnu}_2) < n\widetilde{\mbnu}_2'
\mbSigma_{22}\widetilde{\mbnu}_2-p - 2c \,\}.
\end{eqnarray}
In passing, the condition $n\widetilde{\mbnu}_2'\mbSigma_{22}\widetilde{\mbnu}_2 - p - 2c > 0$
is needed to ensure that $\mathcal{H}^{a}_c$ is not an empty set. Consider
$\mbZ_2 = \mbSigma_{22}\widetilde{\mbnu}_2$; notably, (a) ${\mathcal H}^{a}_c$ is not an empty
set and (b) ${\mathcal H}^{a}_c \cap {\mathcal A}_{U} \neq \varnothing$, which implies that
${\mathcal H}_c \cap {\mathcal A}_{U} \neq \varnothing$. This leads to a contradiction.

\indent Finally, in case (iii) $\widetilde{\mbnu}_1 = {\bf 0}$ and $\widetilde{\mbnu}_2 = {\bf 0}$.
Then, the set ${\mathcal H}^{a}_c$ in equation (2.35)
reduces to
\begin{eqnarray}
{\mathcal H}^{a}_c = \{\, \mbZ \mid n\mbZ_1'\mbSigma_{11:2}^{-1}\mbZ_1
+ n {\mbZ}^{'}_2 {\mbSigma}^{-1}_{22}{\mbZ_2} < -p - 2c \,\}.
\end{eqnarray}
In this case, $ p + 2c < 0$ is required to ensure that $\mathcal{H}^{a}_c$ is not an empty
set. That ${\mathcal H}^{a}_c \cap {\mathcal A}_{U} \neq \varnothing$ can be easily seen,
so ${\mathcal H}_c \cap {\mathcal A}_{U} \neq \varnothing$. This leads to a contradiction.

\indent The discussions of (i), (ii) and (iii) can be taken together to imply that
\begin{align}
\widetilde{\mbnu}_1 \neq {\bf 0} \quad \mbox{and} \quad \widetilde{\mbnu}_2 = {\bf 0}.
\end{align}
Therefore, Rao's $U$-test satisfies the conditions of Theorem 8 of Lehmann ([5], pages 307).
Marden and Perlman [6] have shown that Rao's $U$-test is similar and unbiased, and
the theorem follows by Corollary 2 of Lehmann ([5], page 308).

\vspace{0.3cm}
\def \theequation{3.\arabic{equation}}
\setcounter{equation}{0}
\noindent {\bf 3. Inadmissibility of Hotelling's $T^{2}$-test}
\vspace{0.3cm}

\indent Marden and Perlman ([6], p. 49) pointed out that, ``by utilizing the exponential
structure of the distribution of $(\mXbar, {\bf S})$, the method of Stein [15] and
Schwartz [13] can be applied to reveal the overall $T^{2}$ test is admissible for problem
(1.2). Based on the logarithm of the joint density of $(\mXbar, {\bf S})$, Marden and
Perlman ([6], pages 49-50) claimed that, according to the theorem of Stein [15], the set
(in our notation)
\begin{eqnarray} \nonumber
& \{\, (\mXbar, \mbS) \mid \sup_{(\mbtheta, \mbSigma) \in {\mbTheta}_{2}}
   -\frac{1}{2}n(\mXbar-{\mbtheta})^{'}{\mbSigma}^{-1}(\mXbar-{\mbtheta})
   -\frac{1}{2} (n-1)[\tr {\mbSigma}^{-1}{\bf S}/(n-1)  \\
& \hspace{6cm}
    -\mbox{ln} |{\mbSigma}^{-1}{\bf S}/(n-1)|] \leq  c \,\}
\end{eqnarray}
is an admissible acceptance region in problem (1.2) for any subset
${\mbTheta}_{2} \subset {\mbTheta}_{1}$, where
${\mbTheta}_{1}=\{\, (\mbtheta, \mbSigma) \mid {\mbtheta}_{1}
\ne {\bf 0}, {\mbtheta}_{2}={\bf 0}, {\mbSigma} ~\hbox{p.d.}\}$. Note that the notion in (3.1)
is essentially the same as the one presented in Marden and Perlman ([6], page 50), but omits
the terms for which the parameters and statistics can be separated. The method presented in (3.1)
is easier to handle. However, Marden and Perlman [6] did not offer an analytical proof
for their assertion. Note that the problems considered by Stein [15] and Schwartz [13] are
fully ${\it G}_{0}$-invariant. For the ${\it G}_{0}$-invariant models considered by Stein
[15], ${\mbTheta}_{1}=\{\, (\mbtheta, \mbSigma) \mid {\mbtheta}
\ne {\bf 0}, {\mbSigma} ~\hbox{p.d.}\}$ for the problem of testing
$H^{u}_{0}:{\mbtheta} = {\bf 0}$ against the global alternative
$H^{u}_{1}:{\mbtheta} \ne {\bf 0}$. Note that $\mXbar$ and ${\bf S}$ are independent, and ${\mbmu}$
and ${\mbSigma}$ are orthogonal. Take ${\mbTheta}_{2}=\{\,(\mbtheta, \mbSigma) \mid {\mbtheta}\ne {\bf 0},
{\mbtheta}^{'}{\mbSigma}^{-1}{\mbtheta}=1, {\mbSigma} ~\hbox{p.d.}\}$, substituting the estimator
$(n-1)^{-1}{\bf S}$ of ${\mbSigma}$ into ${\mbSigma}$ and adopting the notation defined in Section 2,
yields
\begin{align} \nonumber
&\{\, (\mXbar, \mbS) \mid \sup_{({\mbtheta} \ne {\bf 0},
    {\mbtheta}^{'}{\mbSigma}^{-1}{\mbtheta}=1, {\mbSigma} ~\hbox {p.d.})}
   -\frac{1}{2}n(\mXbar-{\mbtheta})^{'}{\mbSigma}^{-1}(\mXbar-{\mbtheta})
    -\frac{1}{2} (n-1)[\tr {\mbSigma}^{-1}{\bf S}/(n-1) \\ \nonumber
& \hspace{8cm}
     -\mbox{ln} |{\mbSigma}^{-1}{\bf S}/(n-1)|]   \leq c \,\} \\ \nonumber
& = \{\, (\mXbar, \mbS) \mid \sup_{(\widetilde{\mbnu} \ne {\bf 0},
    \widetilde{\mbnu}^{'}\widetilde{\mbS}^{-1}\widetilde{\mbnu}=1)}
    n(n-1)\widetilde{\mbnu}^{'}\widetilde{\bf X}
    -\frac{1}{2}n(n-1)\widetilde{\bf X}^{'}\widetilde{\mbS}^{-1}\widetilde{\bf X}
    \leq c+\frac{n(p+1)-p}{2} \,\}
\end{align}
\begin{align} \nonumber
& =\{\, ({\mXbar}, {\mbS}) \mid n(n-1)\widetilde{\bf X}^{'}\widetilde{\mbS}^{-1}
   \widetilde{\bf X} \leq k^{'}\,\}  \\  \nonumber
& = \{\, ({\mXbar}, {\mbS}) \mid n(n-1){\mXbar}'{\bf S}^{-1}{\mXbar} \leq k^{'}\,\},
\end{align}
where $\widetilde{\bf X}=({\mXbar}^{'}_{1:2},{\mXbar}^{'}_{2})^{'},
\widetilde{\mbS}=\hbox{diag}({\bf S}_{11.2},{\bf S}_{22})$ and $k^{'}=2c+n(p+1)-p$. The set
(3.2) is equivalent to the acceptance region of Hotelling's $T^{2}$-test.

\indent Note that problem (1.2) is not ${\it Gl}$-invariant, although it is ${\it G}$-invariant.
Marden and Perlman [6] transformed set (3.1) into a ${\it G}$-invariant set to work out set (3.1)
when ${\mbTheta}_{2}={\mbTheta}_{1}$, and reached the conclusion that the ${\it G}$-invariant set
that corresponds to set (3.1) is equivalent to the acceptance region of Hotelling's $T^{2}$-test.
However, in their derivations (Marden and Perlman [6], p. 50) the restriction
${\Delta}_{2}=0$ (${\Delta}_{2}$ defined in Section 1) had to be imposed, thus corresponding to the
assumed condition that ${\mbtheta}_{2}={\bf 0}$ in set (3.1) was overlooked in their new
${\it G}$-invariant set. Rather than focusing only on the ${\it G}$-invariant set, this work directly
determines the form of set (3.1) when ${\mbTheta}_{2}=\{\, (\mbtheta, \mbSigma) \mid {\mbtheta}_{1}
\ne {\bf 0}, {\mbtheta}_{2}={\bf 0}, {\mbtheta}_{1:2}^{'}{\mbSigma}_{11:2}^{-1}{\mbtheta}_{1:2}=
1, {\mbSigma} ~\hbox{p.d.}\}$. Similar to arguments above, (3.1) then becomes
\begin{align}
& \{\, (\mXbar, \mbS) \mid \sup_{(\widetilde{\mbnu}_{1} \ne {\bf 0},
    \widetilde{\mbnu}_{2}={\bf 0},
    \widetilde{\mbnu}_{1}^{'}{\bf S}^{-1}_{11:2}\widetilde{\mbnu}_{1}=1)}
    n(n-1)\widetilde{\mbnu}^{'}\widetilde{\bf X}
   -\frac{1}{2}n(n-1)\widetilde{\bf X}^{'}\widetilde{\mbS}^{-1}\widetilde{\bf X}
  \leq c+\frac{n(p+1)-p}{2} \,\} \\  \nonumber
& = \{\, ({\mXbar}, {\mbS}) \mid
     n(n-1)({\mXbar}'_{1:2}{\bf S}^{-1}_{11:2}{\mXbar}_{1:2}
     -{\mXbar}'_{2}{\bf S}_{22}^{-1}{\mXbar}_{2}) \leq k^{'}\,\}.
\end{align}
Notably, the set (3.3) is also a ${\it G}$-invariant set, but it is not equivalent to the
acceptance region of Hotelling's $T^{2}$-test
\begin{eqnarray}\nonumber
\mathcal {A}_{T^{2}}&=&\{({\mXbar},{\bf S}) \mid
         n(n-1){\mXbar}'{\bf S}^{-1}{\mXbar}\leq k^{*}\}  \\  \nonumber
&=& \{\, ({\mXbar}, {\mbS}) \mid
     n(n-1)({\mXbar}'_{1:2}{\bf S}^{-1}_{11:2}{\mXbar}_{1:2}
     +{\mXbar}'_{2}{\bf S}_{22}^{-1}{\mXbar}_{2}) \leq k^{*}\,\}
\end{eqnarray}
for a suitable $k^{*}$.  Note that, (3.3) is obtained by using the information that
${\mbtheta}_{2}={\bf 0}$, but (3.2) is obtained without using that information.

\indent Due to the fact that $\mXbar$ and ${\bf S}$ are independent, and ${\mbmu}$ and
${\mbSigma}$ are orthogonal; based on the above discussions, an admissible acceptance region
for the problem of testing $H^{u}_{0}:{\mbtheta} = {\bf 0}$ against the global alternative
$H^{u}_{1}:{\mbtheta} \ne {\bf 0}$ can be simply taken as
\begin{align}
& \{\, (\mXbar, \mbS) \mid \sup_{(\widetilde{\mbnu} \ne {\bf 0},
    \widetilde{\mbnu}^{'}\widetilde{\mbS}^{-1}\widetilde{\mbnu}=1)}
    n(n-1)\widetilde{\mbnu}^{'}\widetilde{\bf X} \leq c \,\} \\  \nonumber
& =\{\, ({\mXbar}, {\mbS}) \mid n(n-1)\widetilde{\bf X}^{'}\widetilde{\mbS}^{-1}
   \widetilde{\bf X} \leq c\,\}  \\  \nonumber
& = \{\, ({\mXbar}, {\mbS}) \mid n(n-1){\mXbar}'{\bf S}^{-1}{\mXbar} \leq c\,\},
\end{align}
for a suitable $c$. For problem (1.2), another admissible acceptance region is of the form
\begin{align}
& \{\, (\mXbar, \mbS) \mid \sup_{(\widetilde{\mbnu}_{1} \ne {\bf 0},
    \widetilde{\mbnu}_{2}={\bf 0},
    \widetilde{\mbnu}_{1}^{'}{\bf S}^{-1}_{11:2}\widetilde{\mbnu}_{1}=1)}
    n(n-1)\widetilde{\mbnu}^{'}\widetilde{\bf X}\leq c \,\} \\  \nonumber
& = \{\, ({\mXbar}, {\mbS}) \mid
     n(n-1){\mXbar}'_{1:2}{\bf S}^{-1}_{11:2}{\mXbar}_{1:2} \leq c\,\},
\end{align}
which is the acceptance region of Rao's $U$-test. Instead, based on the unproved assertion (3.1),
we have provided an analytical proof that Rao's $U$-test is admissible for problem (1.2)
by using the Birnbaum-Stein method in Section 2.

\indent Notably, for each ${\bf b} \in R^{p}$, $T^{2}$ and $U$ can be obtained by
maximizing $[n(n-1)]^{1/2}{\bf b}'({\mXbar}^{'}_{1:2},{\mXbar}^{'}_{2})^{'}$ under the
condition that ${\bf b}'\widetilde{\bf S}^{-1}{\bf b}$ is constant over the sets
${\Omega}^{*}_{1}=\{{\bf b}\in R^{p}|~{\bf b}\ne {\bf 0}\}$ and
${\Omega}^{*}_{2}=\{{\bf b}\in R^{p}|~{\bf b}_{1}\ne {\bf 0},~{\bf b}_{2}={\bf 0}\}$,
respectively. Thus, both Hotelling's $T^{2}$-test statistic and Rao's $U$-test
statistic can be constructed by applying the union-intersection (UI) principle of Roy
[12] for the problem of testing $H^{u}_{0}:{\mbtheta} = {\bf 0}$ against the global
alternative $H^{u}_{1}:{\mbtheta} \ne {\bf 0}$, and the problem (1.2) of testing $H_{0}:
{\mbtheta}_{1} = {\bf 0}, {\mbtheta}_{2} = {\bf 0}$ against the alternative
$H_{1}: {\mbtheta}_{1} \ne {\bf 0}, {\mbtheta}_{2}={\bf 0} $, respectively.
Therefore, for problem (1.2) Rao's test based on $U$ may be regarded as a UI test.
Although the Rao $U$-test statistic is constructed by incorporating the information
of ${\mbtheta}_{2}={\bf 0}~(\widetilde{\mbnu}_2={\bf 0})$, the Hotelling $T^{2}$-test
statistic is not thus determined. Therefore, Hotelling's $T^{2}$-test may be reasonably
thought to be dominated by Rao's $U$-test for problem (1.2). This assertion can be
numerically confirmed by the results of Tables 4.1a, 4.1b and 4.1c of Marden and Perlman
[6]. The Birnbaum-Stein method fails to determine whether Hotelling's $T^{2}$-test is
inadmissible for problem (1.2). This shortcoming is overcomed herein by applying Eaton's [3]
basic results to an essentially complete class of test functions for problem (1.2). Let
$\Phi$ be Eaton's essentially complete class of tests, so for any test
${\varphi}^{*}\notin {\Phi}$, there exists a test ${\varphi} \in {\Phi}$ such that
${\varphi}$ is at least as good as ${\varphi}^{*}$.

\vspace{0.3cm}
\noindent {\bf Theorem 3.} {\it For the problem (1.2), Hotelling's $T^{2}$-test is inadmissible}.

\indent {\bf Proof}. Following Eaton [3],  the following is defined.
\begin{eqnarray}
 \Omega_{1}=\{{\mbSigma}^{-1}{\mbtheta}|~{\mbtheta}_{1} \ne {\bf 0},
  {\mbtheta}_{2}={\bf 0}\}\backslash \{{\bf 0}\}.
\end{eqnarray}
Let ${\cal V}\subseteq R^{p}$ be the smallest closed convex cone that contains $\Omega_{1}$.
Then the dual cone of ${\cal V}$ is defined as
\begin{eqnarray}
 {\cal V}^{-}=\{{\bf w}|<{\bf w},{\bf x}> \leq 0,~\forall~ {\bf x}\in {\cal V}\}.
\end{eqnarray}
Notably, that $\Omega_{1}$ is contained in some half-space is not a necessary condition but a
sufficient condition that ensures the dual cone ${\cal V}^{-}$ is a non-empty set. Note that although
${\mbSigma}$ is unknown, but it is fixed. By (3.5) and (3.6), thus we have
\begin{eqnarray}
{\cal V}&=&\{{\mbSigma}^{-1}{\mbtheta}|~{\mbtheta}_{1} \ne {\bf 0},
             {\mbtheta}_{2}={\bf 0}, {\mbSigma}~\mbox {is p.d.}\}  \\ \nonumber
         &=& R^{p_1},
\end{eqnarray}
which is contained in a half-space of $R^{p}$. Similarly, its dual cone is
\begin{eqnarray}
{\cal V}^{-}&=&\{{\bf w}|{\bf x}^{'}{\bf w} \leq 0,~\forall~{\bf x}\in
        {\cal V} \} \\\nonumber
     &=&\{\widetilde{\bf w}|~{\mbtheta}^{'}_{1}{\mbSigma}^{-1}_{11:2}\widetilde{\bf w}_{1}
     \leq  0,~{\mbtheta}_{1} \ne {\bf 0},  {\mbSigma}_{11:2}~\mbox {is p.d.} \},
\end{eqnarray}
where
\begin{eqnarray}
     \widetilde{\bf w}&=&\left[
\begin{array}{cc}
{\bf I}& -\mbSigma_{12}{\mbSigma}^{-1}_{22}\\
{\bf 0}& {\bf I}
\end{array}
\right]{\bf w}  \\ \nonumber
&=&\left(
  \begin{array}{c}
     {\bf w}_{1}-{\mbSigma}_{12}{\mbSigma}^{-1}_{22}{\bf w}_{2}  \\
     {\bf w}_{2}
  \end{array}
  \right)  \\ \nonumber
&\triangleq& \left(
  \begin{array}{c}
     \widetilde{\bf w}_{1}  \\
     \widetilde{\bf w}_{2}
  \end{array}
  \right).
\end{eqnarray}
Therefore, ${\cal V}^{-}=R^{p_2}$, which is an unbounded set.

\indent The acceptance region of Hotelling's $T^{2}$-test is given by
\begin{eqnarray}
\mathcal {A}_{T^{2}}=\{({\mXbar},{\bf S})
   \mid  n(n-1){\mXbar}'{\bf S}^{-1}{\mXbar} \leq {t^{2}_{\alpha}}\},
\end{eqnarray}
where ${t^{2}_{\alpha}}$ is the upper $100\alpha$\% point of the null hypothesis distribution of
$T^{2}$ (which is linked to a F-distribution). For fixed ${\bf S}$, $\mathcal {A}_{T^{2}}$ is an
ellipsoidal set with origin ${\bf 0}$, and is bounded, whereas $\mathcal {V}^{-}$, as shown above,
is still unbounded. Therefore, the proposition 2.1 of Eaton that the dual cone ${\cal V}^{-}$ should
be a subset of the acceptance region of Hotelling's $T^{2}$-test (Eaton [3], section 4, p. 1887)
is not tenable, and thus Hotelling's $T^{2}$-test is not a member of an essentially complete class.

\vspace{0.3cm}
\def \theequation{4.\arabic{equation}}
\setcounter{equation}{0}
\noindent {\bf 4. Whither Rao's $W$-test?}
\vspace{0.3cm}

\indent In passing, both the Hotelling $T^{2}$-test statistic $T^{2}$ and the Rao $W$-test statistic $W$
can be obtained by applying the likelihood ratio principle. The Rao $W$-test statistic is constructed by
incorporating the information that ${\mbtheta}_{2}={\bf 0}$, but the Hotelling $T^{2}$-test statistic is
not thus obtained. For problem (1.10), Marden and Perlman [6] adopted the generalized Bayes approach to
show that Rao's $W$-test is admissible when $0 < \alpha < {\alpha}^{*}$ and is inadmissible when
${\alpha}^{*}  < \alpha < 1$. Section 1 stated that restricting problem (1.2) to ${\it G}$-invariant tests
does not reduce it to problem (1.10), but should reduce it to problem (1.11), and that problem (1.11) and
problem (1.10) differ. Thus, the optimal criteria established for problem (1.10) to draw inferences for
problem (1.2) may lead to conclusions that convey misleading messages. The generalized Bayes approach
of Marden and Perlman [6] can be adopted to characterize in parallel the sufficient and necessary
conditions for the admissibility of the problem (1.11). A situation in which the Rao $W$-test can be
further demonstrated to be a generalized Bayes test, and the corresponding optimality conditions for
the problem (1.11) can be satisfied, can lead to completion of the task. Birnbaum [2], in the context
of complete class type theorems, noted that for testing $H^{u}_{0}:{\mbtheta} ={\bf 0}$ versus
$H^{u}_{1}:{\mbtheta} \ne {\bf 0}$, a test is admissible if and only if it is a generalized Bayes test.
Some admissible tests in the literature are not the Bayes tests for other hypothesis testing
problems (Oosterhoff [8], p.82). For problem (1.2), the set of proper Bayes tests and their weak limits
might only constitute a proper subset of an essentially complete class of tests. On the other hand,
the Birnbaum-Stein method stipulated the convexity assumption for the acceptance regions of tests. However,
for problem (1.2) the acceptance region of Rao's $W$-test is a hyperbolic type set, which is no longer convex.
Therefore, the Birnbaum-Stein method fails to be applicable to Rao's $W$-test. A future study will investigate
the problem of the optimality of Rao's $W$-test.

\vspace{0.3cm}
\def \theequation{5.\arabic{equation}}
\setcounter{equation}{0}
\noindent {\bf 5. Some remarks}
\vspace{0.3cm}

\indent The Hotelling $T^2$-test enjoys many optimal properties of the Neyman-Pearson hypothesis testing theory
when testing against the global alternative. These include similarity, unbiasedness, power monotonicity,
most stringency, uniformly most powerful invariant and alpha-admissibility etc.. However, it is still open
to debate whether the Hotelling $T^{2}$-test is minimax. For the hypothesis testing problem (1.2) we show that
the Hotelling $T^{2}$-test is not a member of an essentially complete class (Eaton [3]), and hence it is no
longer admissible. Moreover, we adopt the Birnbaum-Stein method (Stein [15]) to demonstrate that the Rao
$U$-test is admissible for the hypothesis testing problem (1.2).

\indent Consider the hypotheses
\begin{eqnarray}
H_{0}:{\mbtheta} = {\bf 0}~~\hbox{vs.}~~H_{1}:{\mbtheta} \in {\cal C} \backslash \{{\bf 0}\},
\end{eqnarray}
where ${\cal C}$ denotes a closed convex cone containing a p-dimensional open set. Denote the positive orthant
space by ${\cal O}_{p}^{+} = \{{\mbtheta}\in R^{p}|~{\mbtheta} \geq {\bf 0}\}$. Notice that when ${\cal C}$
is a proper set contained in a halfspace, under a suitable linear transformation the problem in (5.1) can
be reduced to the problem for testing against the positive orthant space with another unknown positive definite
covariance matrix. When ${\cal C}$ is a specific halfspace, then it can be transformed into another halfspace
by a non-singular linear transformation. Hence, without loss of generality it is sufficient to study the cases
in which ${\cal C}$ is the positive orthant space ${\cal O}_{p}^{+}$ and ${\cal C}$ is the halfspace
${\cal H}_{p}^{*}=\{{\mbtheta}\in R^{p}|~{\theta}_{p} \geq 0 \}$. Note that the hypothesis testing problem
(1.2) is a special case of the hypothesis testing problem (5.1). Therefore, we will further study whether
the property of $d$-admissibility of the UIT and LRT for the problem of testing against the closed convex
cone holds in the near future.

\vspace{0.3cm}
\def \theequation{A.\arabic{equation}}
\setcounter{equation}{0}
\noindent {\bf Appendix}
\vspace{0.3cm}

\indent  The following six lemmas are established to prove Theorem 1.

\vspace{0.3cm}
\noindent {\bf Lemma 1}. {\it ${\bf B}^{+}(\mbS)$ is convex on ${\cal S}$}.

\indent {\bf Proof.} For any two given matrices $\mbS$ and $\mbT ~(\in {\cal S})$, let
$\mbPsi^{+}(\alpha) = {\bf B}^{+}(\alpha\mbS + (1 - \alpha)\mbT)$, where $\alpha \in (0, 1)$. Then,
\begin{align*}
\frac{\mathop{d}\, \mbPsi^{+}(\alpha)}{\mathop{d}\, \alpha}
&= - \Bigg\{ (\alpha\mbS + (1 - \alpha)\mbT)^{-1}(\mbS - \mbT)
(\alpha\mbS + (1 - \alpha)\mbT)^{-1} \Bigg.\\
&\qquad \left.
- \left[
\begin{array}{cc}
{\bf 0}& {\bf 0}\\
{\bf 0}& (\alpha\mbS_{22} + (1 - \alpha)\mbT_{22})^{-1}(\mbS_{22} - \mbT_{22})
(\alpha\mbS_{22} + (1 - \alpha)\mbT_{22})^{-1}
\end{array}
\right]
\right\}_.
\end{align*}
Write
\begin{align}\nonumber
\alpha\mbS + (1 - \alpha)\mbT
&= \left[
\begin{array}{cc}
\alpha\mbS_{11} + (1 - \alpha)\mbT_{11}& \alpha\mbS_{12} + (1 - \alpha)\mbT_{12}\\
\alpha\mbS_{21} + (1 - \alpha)\mbT_{21}& \alpha\mbS_{22} + (1 - \alpha)\mbT_{22}
\end{array}
\right]  \\  \nonumber
&\triangleq \left[
\begin{array}{cc}
{\bf E}& {\bf F}\\
{\bf F}'& {\bf G}
\end{array}
\right]
\end{align}
and
\begin{align}  \nonumber
(\alpha\mbS + (1 - \alpha)\mbT)^{-1}
&= {\bf H} + \left[
\begin{array}{cc}
{\bf 0}& {\bf 0}\\
{\bf 0}& {\bf G}^{-1}
\end{array}
\right], \quad \mbox{where }  \\  \nonumber
&{\bf H} = \left(
\begin{array}{c}
{\bf I}\\
-{\bf G}^{-1}{\bf F}'\\
\end{array}
\right)
({\bf E} - {\bf F}{\bf G}^{-1}{\bf F}')^{-1}
\left(
\begin{array}{cc}
{\bf I}& -{\bf F}{\bf G}^{-1}
\end{array}
\right)_.
\end{align}
Since ${\bf H}$ is p.s.d., thus
\[
(\alpha\mbS + (1 - \alpha)\mbT)^{-1} \succeq \left[
\begin{array}{cc}
{\bf 0}& {\bf 0}\\
{\bf 0}& {\bf G}^{-1}
\end{array}
\right]_.
\]
Therefore,
\begin{eqnarray}
\lambda_{max}\left(\left({\bf H} + \left[
\begin{array}{cc}
{\bf 0}& {\bf 0}\\
{\bf 0}& {\bf G}^{-1}
\end{array}
\right]\right)^{-1}
\left[
\begin{array}{cc}
{\bf 0}& {\bf 0}\\
{\bf 0}& {\bf G}^{-1}
\end{array}
\right]
\right) \leq 1,
\end{eqnarray}
where $\lambda_{max}({\bf A})$ denotes the largest eigenvalue of ${\bf A}$. Suppose that
${\bf A} $ is p.d., ${\bf B}$ is p.s.d. and write ${\bf A} = {\bf A}^{1/2}({\bf A}^{1/2})'$,
then note that ${\bf A} \succeq {\bf B}$ implies that
${\bf I}\succeq {\bf A}^{-1/2}{\bf B}({\bf A}^{-1/2})'$.
Thus, $1 =\lambda_{max}({\bf I}) \geq \lambda_{max}({\bf A}^{-1/2} {\bf B}({\bf A}^{-1/2})')
=\lambda_{max}({\bf B}{\bf A}^{-1}) = \lambda_{max}({\bf A}^{-1} {\bf B})$. Note that,

\begin{align*}
\frac{\mathop{d}^2\ \mbPsi^{+}(\alpha)}{\mathop{d}\ \alpha^2}
&= 2 \Bigg\{
(\alpha\mbS + (1 - \alpha)\mbT)^{-1}(\mbS - \mbT)(\alpha\mbS + (1 - \alpha)\mbT)^{-1}
(\mbS - \mbT)(\alpha\mbS + (1 - \alpha)\mbT)^{-1} \Bigg.\\
&\quad \Bigg. - \left[
\begin{array}{cc}
{\bf 0}& {\bf 0}\\
{\bf 0}& {\scriptstyle (\alpha\mbS_{22}
+ (1 - \alpha)\mbT_{22})^{-1}(\mbS_{22} - \mbT_{22})
(\alpha\mbS_{22} + (1 - \alpha)\mbT_{22})^{-1}(\mbS_{22}
- \mbT_{22})(\alpha\mbS_{22} + (1 - \alpha)\mbT_{22})^{-1}}
\end{array}
\right] \Bigg\}  \\
&= 2 \Bigg\{(\alpha\mbS + (1 - \alpha)\mbT)^{-1}(\mbS - \mbT)(\alpha\mbS
+ (1 - \alpha)\mbT)^{-1}(\mbS - \mbT)(\alpha\mbS + (1 - \alpha)\mbT)^{-1} \Bigg.\\
&\qquad - \Bigg. \left[
\begin{array}{cc}
{\bf 0}& {\bf 0}\\
{\bf 0}& {\bf G}^{-1}
\end{array}
\right]
(\mbS - \mbT) \left[
\begin{array}{cc}
{\bf 0}& {\bf 0}\\
{\bf 0}& {\bf G}^{-1}
\end{array}
\right]
(\mbS - \mbT) \left[
\begin{array}{cc}
{\bf 0}& {\bf 0}\\
{\bf 0}& {\bf G}^{-1}
\end{array}
\right]
\Bigg\}\\
&\succeq 2 \Bigg\{ (\alpha\mbS + (1 - \alpha)\mbT)^{-1}(\mbS - \mbT)
(\alpha\mbS + (1 - \alpha)\mbT)^{-1}(\mbS - \mbT)(\alpha\mbS + (1 - \alpha)\mbT)^{-1}\\
&\quad - \Bigg. \left[
\begin{array}{cc}
{\bf 0}& {\bf 0}\\
{\bf 0}& {\bf G}^{-1}
\end{array}
\right]
(\mbS - \mbT)(\alpha\mbS - (1 - \alpha)\mbT)^{-1}(\mbS - \mbT) \left[
\begin{array}{cc}
{\bf 0}& {\bf 0}\\
{\bf 0}& {\bf G}^{-1}
\end{array}
\right]
\Bigg\}\\
&= 2 \Bigg\{\, (\alpha\mbS + (1 - \alpha)\mbT)^{-1}{\bf K}(\alpha\mbS
+ (1 - \alpha)\mbT)^{-1} - \left[
\begin{array}{cc}
{\bf 0}& {\bf 0}\\
{\bf 0}& {\bf G}^{-1}
\end{array}
\right]
{\bf K} \left[
\begin{array}{cc}
{\bf 0}& {\bf 0}\\
{\bf 0}& {\bf G}^{-1}
\end{array}
\right] \,\Bigg\},
\end{align*}
where
\[
{\bf K} = (\mbS - \mbT)(\alpha\mbS + (1 - \alpha)\mbT)^{-1}(\mbS - \mbT).
\]

\indent Next, compare the matrices
\[
(\alpha\mbS + (1 - \alpha)\mbT)^{-1}{\bf K}(\alpha\mbS + (1 - \alpha)\mbT)^{-1}
\]
and
\[
\left[
\begin{array}{cc}
{\bf 0}& {\bf 0}\\
{\bf 0}& {\bf G}^{-1}
\end{array}
\right] {\bf K} \left[
\begin{array}{cc}
{\bf 0}& {\bf 0}\\
{\bf 0}& {\bf G}^{-1}
\end{array}
\right]_,
\]
that is,
\[
\left(
{\bf H} + \left[
\begin{array}{cc}
{\bf 0}& {\bf 0}\\
{\bf 0}& {\bf G}^{-1}
\end{array}
\right] \right)
{\bf K} \left( {\bf H}
+ \left[
\begin{array}{cc}
{\bf 0}& {\bf 0}\\
{\bf 0}& {\bf G}^{-1}
\end{array}
\right]
\right)
\]
and
\[
\left[
\begin{array}{cc}
{\bf 0}& {\bf 0}\\
{\bf 0}& {\bf G}^{-1}
\end{array} \right] {\bf K} \left[
\begin{array}{cc}
{\bf 0}& {\bf 0}\\
{\bf 0}& {\bf G}^{-1}
\end{array}
\right]_.
\]
Consider the new matrix ${\bf L}{\bf L}'$, where
\[
{\bf L} = {\bf K}^{-1/2} \left( {\bf H} + \left[
\begin{array}{cc}
{\bf 0}& {\bf 0}\\
{\bf 0}& {\bf G}^{-1}
\end{array}
\right] \right)^{-1} \left[
\begin{array}{cc}
{\bf 0}& {\bf 0}\\
{\bf 0}& {\bf G}^{-1}
\end{array}
\right] {\bf K}^{1/2}.
\]
The matrix ${\bf L}$ can be rewritten as ${\bf L} = {\bf P}\mbLambda{\bf Q}'$, where $\mbLambda$
denotes the diagonal matrix of eigenvalues of ${\bf L}$ and ${\bf P}$, ${\bf Q} \in \mathcal{Q}(p)$,
the group of $p \times p$ orthogonal matrices. Thus, ${\bf L}{\bf L}' = {\bf P}\mbLambda^2{\bf P}'$.
Therefore, $\lambda_{max}({\bf L}{\bf L}') = \lambda_{max}({\bf P}\mbLambda^2{\bf P}')
= \lambda_{max}(\mbLambda^2{\bf P}'{\bf P}) = \lambda_{max}(\mbLambda)^2
= [\lambda_{max}(\mbLambda)]^2 = [\lambda_{max}({\bf L})]^2$. Notably,
\begin{align*}
\lambda_{max}({\bf L})
&= \lambda_{max} \left({\bf K}^{-1/2} \left( {\bf H} + \left[
\begin{array}{cc}
{\bf 0}& {\bf 0}\\
{\bf 0}& {\bf G}^{-1}
\end{array}
\right] \right)^{-1} \left[
\begin{array}{cc}
{\bf 0}& {\bf 0}\\
{\bf 0}& {\bf G}^{-1}
\end{array}
\right] {\bf K}^{1/2} \right)\\
&= \lambda_{max} \left( \left({\bf H} + \left[
\begin{array}{cc}
{\bf 0}& {\bf 0}\\
{\bf 0}& {\bf G}^{-1}
\end{array}
\right] \right)^{-1} \left[
\begin{array}{cc}
{\bf 0}& {\bf 0}\\
{\bf 0}& {\bf G}^{-1}
\end{array}
\right] \right) \\ \nonumber
& \leq 1,
\end{align*}
the last inequality follows from the inequality (A.1). Thus,
$\lambda_{max}({\bf L}{\bf L}') \leq 1$. Furthermore, notice that
\begin{align*}
\lambda_{max}({\bf L}{\bf L}')
&= \lambda_{max} \Bigg( \big[ (\alpha\mbS + (1 - \alpha)\mbT)^{-1}(\mbS - \mbT)(\alpha\mbS
+ (1 - \alpha)\mbT)^{-1}(\mbS - \mbT)(\alpha\mbS + (1 - \alpha)\mbT)^{-1} \big]^{-1}\\
&\qquad \left[
\begin{array}{cc}
{\bf 0}& {\bf 0}\\
{\bf 0}& {\bf G}^{-1}
\end{array}
\right]
(\mbS - \mbT)(\alpha\mbS + (1 - \alpha)\mbT)^{-1}(\mbS - \mbT) \left[
\begin{array}{cc}
{\bf 0}& {\bf 0}\\
{\bf 0}& {\bf G}^{-1}
\end{array}
\right] \Bigg)_.
\end{align*}
Thus,
\begin{align*}
&(\alpha\mbS + (1 - \alpha)\mbT)^{-1}(\mbS - \mbT)(\alpha\mbS + (1 - \alpha)\mbT)^{-1}
 (\mbS - \mbT)(\alpha\mbS + (1 - \alpha)\mbT)^{-1}\\
&\qquad \qquad \succeq \left[
\begin{array}{cc}
{\bf 0}& {\bf 0}\\
{\bf 0}& {\bf G}^{-1}
\end{array}
\right]
(\mbS - \mbT)(\alpha\mbS + (1 - \alpha)\mbT)^{-1}(\mbS - \mbT) \left[
\begin{array}{cc}
{\bf 0}& {\bf 0}\\
{\bf 0}& {\bf G}^{-1}
\end{array}
\right]_.
\end{align*}
Therefore,
\[
\frac{\mathop{d}^2 \mbPsi^{+}(\alpha)}{\mathop{d} \alpha^2} \succeq {\bf 0}
\]
and hence
\[
{\bf B}^{+}(\mbS) \mbox{ is convex on }  {\cal S}.
\]

\vspace{0.3cm}
\noindent {\bf Lemma 2}. {\it ${\bf B}(\mbS)$ is the Moore-Penrose generalized inverse
of ${\bf B}^{+}(\mbS)$}.

\indent {\bf Proof.} It can be easily shown that ${\bf B}(\mbS)$ satisfies the following conditions:
\begin{enumerate}
\renewcommand{\labelenumi}{(\roman{enumi})}
\item ${\bf B}(\mbS){\bf B}^{+}(\mbS){\bf B}(\mbS) = {\bf B}(\mbS)$,
\item ${\bf B}^{+}(\mbS){\bf B}(\mbS){\bf B}^{+}(\mbS) = {\bf B}^{+}(\mbS)$,
\item ${\bf B}(\mbS){\bf B}^{+}(\mbS) = \left( {\bf B}(\mbS){\bf B}^{+}(\mbS) \right)'$,
\item ${\bf B}^{+}(\mbS){\bf B}(\mbS) = \big( {\bf B}^{+}(\mbS){\bf B}(\mbS) \big)'.$
\end{enumerate}

\vspace{0.3cm}
\noindent {\bf Lemma 3}. {\it For any two given matrices $\mbS$ and $\mbT ~(\in {\cal S})$, let
$\mbPsi^{+}(\alpha) = {\bf B}^{+}(\alpha\mbS + (1 - \alpha)\mbT)$, where $\alpha \in (0, 1)$.
Then $\frac{\mathop{d}^2 \mbPsi^{+}(\alpha)}{\mathop{d} \alpha^2} \succeq \left.
\frac{\mathop{d}^2 \mbPsi^{+}(\alpha)}{\mathop{d} \alpha^2} \right|_{\alpha= \alpha_0}$,
$\forall\, \alpha \in (0, 1)$, where $\alpha_0$ is the stationary point of $\mbPsi^{+}(\alpha)$.
Let $\mbPsi(\alpha)$ be the Moore-Penrose generalized inverse of $\mbPsi^{+}(\alpha)$,
then $\frac{\mathop{d}^{2} \mbPsi(\alpha)}{\mathop{d} \alpha^2}
\preceq
\left. \frac{\mathop{d}^2 \mbPsi(\alpha)}{\mathop{d} \alpha^2} \right|_{\alpha
= \alpha_0}$, $\forall\, \alpha \in (0, 1)$}.

\indent {\bf Proof.}
It is easily to see that $\mbPsi^{+}(\alpha)$ is continuous and hence differentiable,
$\forall\, \alpha \in (0, 1)$. First note that $\mbPsi^{+}(\alpha)$ is neither a linear
function nor a quadratic function of $\alpha$. From the proof of Lemma 1, we note that
$\mbPsi^{+}(\alpha)$ is not a monotone function of $\alpha$. By Lemma 1, $\mbPsi^{+}(\alpha)$
is convex in $\alpha$. Thus, the stationary point of $\mbPsi^{+}(\alpha)$ exists and unique.
Suppose $\alpha_0$ be the stationary (critical) point, then $\mbPsi^{+}(\alpha) \succeq
\mbPsi^{+}(\alpha_0) \succeq {\bf 0}$, $\forall\, \alpha \in (0, 1)$. This implies that
$\mbPsi^{+}(\alpha) \succeq (1 - \alpha)^2\mbPsi^{+}(\alpha_0) = (\alpha^2 - 2\alpha +
1)\mbPsi^{+}(\alpha_0)$, $\forall\, \alpha \in (0, 1)$. Note that $\mbPsi^{+}(\alpha_0)$ is p.s.d.,
thus there exists a quadratic convex function ${\mbUpsilon}^{+}(\alpha)$ such that
$\mbPsi^{+}(\alpha) \succeq {\mbUpsilon}^{+}(\alpha)$, $\forall\, \alpha \in (0, 1)$.
Hence for sufficiently small $h$, there exists a quadratic convex function
${\mbUpsilon}_0^{+}(\alpha)$ such that
$\mbPsi^{+}(\alpha)\succeq {\mbUpsilon}_0^{+}(\alpha)$, $\forall\, \alpha \in (0, 1)$, and
$\mbPsi^{+}(\alpha_0) = {\mbUpsilon}_0^{+}(\alpha_0)$, $\mbPsi^{+}(\alpha_0 + h)
= {\mbUpsilon}_0^{+}(\alpha_0 + h)$ and
$\mbPsi^{+}(\alpha_0 + 2h) = {\mbUpsilon}_0^{+}(\alpha_0 + 2h)$. Notably,
\begin{align*}
\frac{\mathop{\mathrm{d}}^2 \mbPsi^{+}(\alpha)}{\mathop{\mathrm{d}} \alpha^2}- \left.
\frac{\mathop{\mathrm{d}}^2 \mbPsi^{+}(\alpha)}{\mathop{\mathrm{d}} \alpha^2}
\right|_{\alpha = \alpha_0}
&= \lim_{h \to 0} \frac{\mbPsi^{+}(\alpha + 2h) - 2\mbPsi^{+}(\alpha + h) +
\mbPsi^{+}(\alpha)}{h^2}\\ &\qquad - \lim_{h \to 0} \frac{\mbPsi^{+}(\alpha_0 + 2h) -
2\mbPsi^{+}(\alpha_0 + h) + \mbPsi^{+}(\alpha_0)}{h^2}   \\
&= \lim_{h \to 0}
\frac{\mbPsi^{+}(\alpha + 2h) - 2\mbPsi^{+}(\alpha + h) + \mbPsi^{+}(\alpha)}{h^2}\\
&\qquad -\lim_{h \to 0}\frac{{\mbUpsilon}_0^{+}(\alpha_0 + 2h) -
2{\mbUpsilon}_0^{+}(\alpha_0 + h) + {\mbUpsilon}_0^{+}(\alpha_0)}{h^2}.
\end{align*}
For any fixed $\alpha \in (0, 1)$, there exists a quadratic convex function
${\mbUpsilon}_{\star}^{+}(\alpha)$ such that ${\mbUpsilon}_{\star}^{+}(\alpha) \succeq
{\mbUpsilon}_0^{+}(\alpha)$, and ${\mbUpsilon}_{\star}^{+}(\alpha) = \mbPsi^{+}(\alpha)$,
${\mbUpsilon}_{\star}^{+}(\alpha + h) = \mbPsi^{+}(\alpha + h)$ and
${\mbUpsilon}_{\star}^{+}(\alpha + 2h) = \mbPsi^{+}(\alpha + 2h)$ for arbitrary small $h$.
Thus,
\begin{align*}
\frac{\mathop{\mathrm{d}}^2 \mbPsi_{}^{+}(\alpha)}{\mathop{\mathrm{d}} \alpha^2}
- \left. \frac{\mathop{\mathrm{d}}^2 \mbPsi^{+}(\alpha)}{\mathop{\mathrm{d}} \alpha^2}
\right|_{\alpha = \alpha_0}
&= \lim_{h \to 0} \frac{{\mbUpsilon}_{\star}^{+}(\alpha + 2h) - 2{\mbUpsilon}_{\star}^{+}
(\alpha + h)+ {\mbUpsilon}_{\star}^{+}(\alpha)}{h^2}\\
&\qquad - \lim_{h \to 0} \frac{{\mbUpsilon}_0^{+}(\alpha_0 + 2h) - 2{\mbUpsilon}_0^{+}
(\alpha_0 + h)+ {\mbUpsilon}_0^{+}(\alpha_0)}{h^2}\\
&  = \frac{\mathop{\mathrm{d}}^2 {\mbUpsilon}_{\star}^{+}
   (\alpha)}{\mathop{\mathrm{d}}\alpha^2}- \left.
\frac{\mathop{\mathrm{d}}^2 {\mbUpsilon}_0^{+}(\alpha)}{\mathop{\mathrm{d}} \alpha^2}
\right|_{\alpha = \alpha_0}\\
& \succeq \frac{\mathop{\mathrm{d}}^2 {\mbUpsilon}_0^{+}
    (\alpha)}{\mathop{\mathrm{d}}\alpha^2}- \left. \frac{\mathop{\mathrm{d}}^2
{\mbUpsilon}_0^{+}(\alpha)}{\mathop{\mathrm{d}} \alpha^2}
 \right|_{\alpha = \alpha_0}= {\bf 0}, \quad \forall\, \alpha \in (0, 1).
\end{align*}
Therefore,
\[
\frac{\mathop{\mathrm{d}}^2 \mbPsi^{+}(\alpha)}{\mathop{\mathrm{d}} \alpha^2}
\succeq \left. \frac{\mathop{\mathrm{d}}^2 \mbPsi^{+}(\alpha)}{\mathop{\mathrm{d}}
\alpha^2} \right|_{\alpha = \alpha_0}, \quad \forall\, \alpha \in (0, 1).
\]
Similarly, $\mbPsi(\alpha)$ is the Moore-Penrose generalized inverse of $\mbPsi^{+}(\alpha)$,
and so $\mbPsi(\alpha_0) \succeq \mbPsi(\alpha) \succeq {\bf 0}$, $\forall\, \alpha \in (0, 1)$.
This implies that $(1 + 2\alpha - \alpha^2)\mbPsi(\alpha_0) \succeq \mbPsi(\alpha)$,
$\forall\, \alpha \in (0, 1)$. Thus, there exists a quadratic concave function ${\mbUpsilon}(\alpha)$
such that ${\mbUpsilon}(\alpha)\succeq \mbPsi(\alpha)$, $\forall\, \alpha \in (0, 1)$.
Parallel arguments as in the case $\mbPsi^{+}(\alpha)$, we may also conclude that
\[
\frac{\mathop{\mathrm{d}}^2 \mbPsi(\alpha)}{\mathop{\mathrm{d}} \alpha^2} \preceq \left.
\frac{\mathop{\mathrm{d}}^2 \mbPsi(\alpha)}{\mathop{\mathrm{d}} {\alpha}^2}
\right|_{\alpha = \alpha_0}, \forall\, \alpha \in (0, 1).
\]

\vspace{0.3cm}
\noindent {\bf Lemma 4}. {\it ${\bf B}(\mbS)$ is concave on ${\cal S}$}.

\indent {\bf Proof.} Rewrite ${\bf B}^{+}(\mbS)$ in (2.2) as
\[
{\bf B}^{+}(\mbS) = \left(
\begin{array}{c}
{\bf I}\\
-\mbS_{22}^{-1}\mbS_{21}
\end{array}
\right) \mbS_{11:2}^{-1} \left(
\begin{array}{cc}
{\bf I}& -\mbS_{12}\mbS_{22}^{-1}
\end{array}
\right)
= {\bf C}(\mbS){\bf D}^{-1}(\mbS){\bf C}'(\mbS),
\]
where
\[
{\bf C}(\mbS) = \left(
\begin{array}{c}
{\bf I}\\
-\mbS_{22}^{-1}\mbS_{21}
\end{array}
\right)
\mbox{ and } {\bf D}(\mbS) = \mbS_{11:2}.
\]
Note that\\
{\Large\textcircled{\raisebox{-.3pt}{\normalsize 1}}}
${\bf B}(\mbS)
= {\bf C}(\mbS) \left({\bf C}'(\mbS){\bf C}(\mbS) \right)^{-1}{\bf D}(\mbS)
\left({\bf C}'(\mbS){\bf C}(\mbS) \right)^{-1}{\bf C}'(\mbS)$\\
{\Large\textcircled{\raisebox{-.3pt}{\normalsize 2}}}
${\bf B}^{+}(\mbS){\bf B}(\mbS) = {\bf B}(\mbS){\bf B}^{+}(\mbS)
={\bf C}(\mbS)({\bf C}'(\mbS){\bf C}(\mbS))^{-1}{\bf C}'(\mbS).$\\
Let $\mbPsi^{+}(\alpha) = {\bf B}^{+}(\alpha\mbS + (1 - \alpha)\mbT)
= {\bf M}{\bf N}^{-1}{\bf M}$, where ${\bf M} = {\bf C}(\alpha\mbS
+ (1 - \alpha)\mbT)$, ${\bf N} = {\bf D}(\alpha\mbS + (1 - \alpha)\mbT)$ and
$\alpha \in (0,1)$. Then, its Moore-Penrose generalized inverse is of the form
\[
\mbPsi(\alpha) = {\bf B}(\alpha\mbS + (1 - \alpha)\mbT)
= {\bf M}({\bf M}'{\bf M})^{-1}{\bf N}({\bf M}'{\bf M})^{-1}{\bf M}'.
\]
Notably,
\[
{\bf B}^{+}(\mbS){\bf B}(\mbS){\bf B}^{+}(\mbS) = {\bf B}^{+}(\mbS) {\implies}
\mbPsi^{+}(\alpha)\mbPsi(\alpha)\mbPsi^{+}(\alpha) = \mbPsi^{+}(\alpha).
\]
Thus,
\begin{align}
\frac{\mathop{\mathrm{d}} \mbPsi^{+}(\alpha)}{\mathop{\mathrm{d}} \alpha}\mbPsi(\alpha)
\mbPsi^{+}(\alpha)+ \mbPsi^{+}(\alpha)\frac{\mathop{\mathrm{d}} \mbPsi(\alpha)}
{\mathop{\mathrm{d}}\alpha}\mbPsi^{+}(\alpha) + \mbPsi^{+}(\alpha)\mbPsi(\alpha)
\frac{\mathop{\mathrm{d}}\mbPsi^{+}(\alpha)}{\mathop{\mathrm{d}} \alpha}
= \frac{\mathop{\mathrm{d}}\mbPsi^{+}(\alpha)}{\mathop{\mathrm{d}} \alpha}
\end{align}
and
\begin{align}
\frac{\mathop{\mathrm{d}}^2 \mbPsi^{+}(\alpha)}{\mathop{\mathrm{d}} \alpha}
& = \frac{\mathop{\mathrm{d}}^2 \mbPsi^{+}(\alpha)}{\mathop{\mathrm{d}}
\alpha}\mbPsi(\alpha)\mbPsi^{+}(\alpha) + \frac{\mathop{\mathrm{d}}
\mbPsi^{+}(\alpha)}{\mathop{\mathrm{d}}\alpha}
\frac{\mathop{\mathrm{d}} \mbPsi(\alpha)}{\mathop{\mathrm{d}} \alpha}\mbPsi^{+}(\alpha)
+ \frac{\mathop{\mathrm{d}} \mbPsi^{+}(\alpha)}{\mathop{\mathrm{d}} \alpha}\mbPsi(\alpha)
\frac{\mathop{\mathrm{d}} \mbPsi^{+}(\alpha)}{\mathop{\mathrm{d}} \alpha}  \\ \nonumber
&\qquad + \frac{\mathop{\mathrm{d}} \mbPsi^{+}(\alpha)}{\mathop{\mathrm{d}} \alpha}
\frac{\mathop{\mathrm{d}} \mbPsi(\alpha)}{\mathop{\mathrm{d}} \alpha}\mbPsi^{+}(\alpha)
+ \mbPsi^{+}(\alpha)\frac{\mathop{\mathrm{d}}^2 \mbPsi(\alpha)}{\mathop{\mathrm{d}}
\alpha^2}\mbPsi^{+}(\alpha)   \\  \nonumber
&\qquad + \mbPsi^{+}(\alpha)\frac{\mathop{\mathrm{d}}
\mbPsi(\alpha)}{\mathop{\mathrm{d}} \alpha}
\frac{\mathop{\mathrm{d}} \mbPsi^{+}(\alpha)}{\mathop{\mathrm{d}} \alpha}
+ \frac{\mathop{\mathrm{d}} \mbPsi^{+}(\alpha)}{\mathop{\mathrm{d}} \alpha}\mbPsi(\alpha)
\frac{\mathop{\mathrm{d}} \mbPsi^{+}(\alpha)}{\mathop{\mathrm{d}} \alpha} \\  \nonumber
&\qquad + \mbPsi^{+}(\alpha)\frac{\mathop{\mathrm{d}} \mbPsi(\alpha)}{\mathop{\mathrm{d}}
 \alpha}\frac{\mathop{\mathrm{d}} \mbPsi^{+}(\alpha)}{\mathop{\mathrm{d}} \alpha}
+ \mbPsi^{+}(\alpha)\mbPsi(\alpha)\frac{\mathop{\mathrm{d}}^2 \mbPsi^{+}(\alpha)}
{\mathop{\mathrm{d}}\alpha^2}.
\end{align}
By the results of (A.2) and (A.3), then
\begin{align*}
\mbPsi^{+}(\alpha)\frac{\mathop{\mathrm{d}}^2 \mbPsi(\alpha)}{\mathop{\mathrm{d}}
\alpha^2}\mbPsi^{+}(\alpha) &= \frac{\mathop{\mathrm{d}}^2 \mbPsi(\alpha)}
{\mathop{\mathrm{d}}\alpha^2} -\mbPsi^{+}(\alpha)\mbPsi(\alpha)
\frac{\mathop{\mathrm{d}}^2 \mbPsi^{+}
(\alpha)}{\mathop{\mathrm{d}}\alpha^2} - \frac{\mathop{\mathrm{d}}^2 \mbPsi^{+}
(\alpha)}{\mathop{\mathrm{d}}\alpha^2}\mbPsi(\alpha)\mbPsi^{+}(\alpha)\\ &\quad -
2\frac{\mathop{\mathrm{d}}\mbPsi^{+}(\alpha)}{\mathop{\mathrm{d}} \alpha}
\frac{\mathop{\mathrm{d}} \mbPsi(\alpha)}{\mathop{\mathrm{d}} \alpha}\mbPsi^{+}(\alpha)
- 2\frac{\mathop{\mathrm{d}} \mbPsi^{+}(\alpha)}{\mathop{\mathrm{d}} \alpha}\mbPsi(\alpha)
\frac{\mathop{\mathrm{d}} \mbPsi^{+}(\alpha)}{\mathop{\mathrm{d}} \alpha}\\
&\qquad - 2\mbPsi^{+}(\alpha)\frac{\mathop{\mathrm{d}} \mbPsi(\alpha)}{\mathop{\mathrm{d}}
\alpha} \frac{\mathop{\mathrm{d}} \mbPsi^{+}(\alpha)}{\mathop{\mathrm{d}} \alpha}.
\end{align*}
Since $\mbPsi^{+}(\alpha)\mbPsi(\alpha) = {\bf M}({\bf M}'{\bf M})^{-1}{\bf M}'$,
thus
\[
\mbPsi^{+}(\alpha)\frac{\mathop{\mathrm{d}} \mbPsi(\alpha)}{\mathop{\mathrm{d}} \alpha}
= -\frac{\mathop{\mathrm{d}} \mbPsi^{+}(\alpha)}{\mathop{\mathrm{d}} \alpha}
\mbPsi(\alpha) + \frac{\mathop{\mathrm{d}}}{\mathop{\mathrm{d}} \alpha}
[{\bf M}({\bf M}'{\bf M})^{-1}{\bf M}'].
\]
Therefore,
\begin{align}
\mbPsi^{+}(\alpha)\frac{\mathop{\mathrm{d}}^2 \mbPsi(\alpha)}{\mathop{\mathrm{d}}
\alpha^2}\mbPsi^{+}(\alpha) &= \frac{\mathop{\mathrm{d}}^2 \mbPsi^{+}(\alpha)}
{\mathop{\mathrm{d}}\alpha^2} - \mbPsi^{+}(\alpha)\mbPsi(\alpha)
\frac{\mathop{\mathrm{d}}^2 \mbPsi^{+}(\alpha)}{\mathop{\mathrm{d}}
\alpha^2} - \frac{\mathop{\mathrm{d}}^2 \mbPsi^{+}(\alpha)}{\mathop{\mathrm{d}}
\alpha^2}\mbPsi(\alpha)\mbPsi^{+}(\alpha)  \\ \nonumber
&\quad + 2\frac{\mathop{\mathrm{d}}
\mbPsi^{+}(\alpha)}{\mathop{\mathrm{d}} \alpha}\mbPsi(\alpha)
\frac{\mathop{\mathrm{d}} \mbPsi^{+}(\alpha)}{\mathop{\mathrm{d}} \alpha}
- 2\frac{\mathop{\mathrm{d}}}{\mathop{\mathrm{d}} \alpha}[{\bf M}({\bf M}'{\bf M})^{-1}
{\bf M}']\frac{\mathop{\mathrm{d}} \mbPsi^{+}(\alpha)}{\mathop{\mathrm{d}} \alpha} \\
\nonumber
&\qquad - 2\frac{\mathop{\mathrm{d}} \mbPsi^{+}(\alpha)}{\mathop{\mathrm{d}} \alpha}
\frac{\mathop{\mathrm{d}}}{\mathop{\mathrm{d}} \alpha}[{\bf M}({\bf M}'{\bf M})^{-1}
{\bf M}'].
\end{align}
Notably,
\begin{enumerate}
\renewcommand{\labelenumi}{(\roman{enumi})}

\item
\begin{align*}
\frac{\mathop{\mathrm{d}}^2 \mbPsi^{+}(\alpha)}{\mathop{\mathrm{d}} \alpha^2}
&= \frac{\mathop{\mathrm{d}}^2 {\bf M}}{\mathop{\mathrm{d}} \alpha^2}{\bf N}^{-1}{\bf M}'
+ {\bf M}{\bf N}^{-1}\frac{\mathop{\mathrm{d}}^2 {\bf M}'}{\mathop{\mathrm{d}} \alpha^2}
+ 2\frac{\mathop{\mathrm{d}} {\bf M}}{\mathop{\mathrm{d}} \phi}{\bf N}^{-1}
\frac{\mathop{\mathrm{d}} {\bf M}'}{\mathop{\mathrm{d}} \alpha} \\
&- 2\frac{\mathop{\mathrm{d}} {\bf M}}{\mathop{\mathrm{d}} \alpha}{\bf N}^{-1}
\frac{\mathop{\mathrm{d}} {\bf M}}{\mathop{\mathrm{d}} \alpha}{\bf N}^{-1}{\bf M}'
 - 2{\bf M}{\bf N}^{-1}\frac{\mathop{\mathrm{d}} {\bf N}}{\mathop{\mathrm{d}}
\alpha}{\bf N}^{-1}\frac{\mathop{\mathrm{d}} {\bf M}'}{\mathop{\mathrm{d}} \alpha}\\
&+ 2{\bf M}{\bf N}^{-1}\frac{\mathop{\mathrm{d}} {\bf N}}{\mathop{\mathrm{d}} \alpha}
{\bf N}^{-1}\frac{\mathop{\mathrm{d}} {\bf N}}{\mathop{\mathrm{d}} \alpha}{\bf N}^{-1}
{\bf M}' -{\bf M}{\bf N}^{-1}\frac{\mathop{\mathrm{d}}^2 {\bf N}}{\mathop{\mathrm{d}}
\alpha^2} {\bf N}^{-1}{\bf M}'\\
&= \frac{\mathop{\mathrm{d}}^2 {\bf M}}{\mathop{\mathrm{d}} \alpha^2}{\bf N}^{-1}{\bf M}'
+ {\bf M}{\bf N}^{-1}\frac{\mathop{\mathrm{d}}^2 {\bf M}'}{\mathop{\mathrm{d}} \alpha^2}
- {\bf M}{\bf N}^{-1}\frac{\mathop{\mathrm{d}}^2 {\bf N}}{\mathop{\mathrm{d}} \alpha^2}
{\bf N}^{-1}{\bf M}'\\
&\quad -2({\bf M}{\bf N}^{-1}\frac{\mathop{\mathrm{d}} {\bf N}}{\mathop{\mathrm{d}} \alpha}
- \frac{\mathop{\mathrm{d}} {\bf M}}{\mathop{\mathrm{d}} \alpha}){\bf N}^{-1}
({\bf M}{\bf N}^{-1}\frac{\mathop{\mathrm{d}} {\bf N}}{\mathop{\mathrm{d}} \alpha}
- \frac{\mathop{\mathrm{d}} {\bf M}}{\mathop{\mathrm{d}} \alpha})',
\end{align*}

\item
\begin{align*}
\mbPsi^{+}(\alpha)\mbPsi(\alpha)\frac{\mathop{\mathrm{d}}^2 \mbPsi^{+}(\alpha)}
{\mathop{\mathrm{d}} \alpha^2}
&= {\bf M}({\bf M}'{\bf M})^{-1}\frac{\mathop{\mathrm{d}}^2 {\bf M}}{\mathop{\mathrm{d}}
 \alpha^2}{\bf N}^{-1}{\bf M}'
 +{\bf M}{\bf N}^{-1}\frac{\mathop{\mathrm{d}}^2 {\bf M}'}{\mathop{\mathrm{d}} \alpha^2}\\
&+ 2{\bf M}({\bf M}'{\bf M})^{-1}{\bf M}'\frac{\mathop{\mathrm{d}} {\bf M}}
{\mathop{\mathrm{d}} \alpha}{\bf N}^{-1}\frac{\mathop{\mathrm{d}}
{\bf M}'}{\mathop{\mathrm{d}} \alpha}\\ & -2{\bf M}({\bf M}'{\bf M})^{-1}{\bf
M}'\frac{\mathop{\mathrm{d}} {\bf M}}{\mathop{\mathrm{d}} \alpha}{\bf N}^{-1}
\frac{\mathop{\mathrm{d}} {\bf N}}{\mathop{\mathrm{d}} \alpha}{\bf N}^{-1}{\bf M}'
 -2{\bf M}{\bf N}^{-1}\frac{\mathop{\mathrm{d}} {\bf N}}{\mathop{\mathrm{d}}
\alpha}{\bf N}^{-1}\frac{\mathop{\mathrm{d}} {\bf M}'}{\mathop{\mathrm{d}} \alpha}\\
&+ 2{\bf M}{\bf N}^{-1}\frac{\mathop{\mathrm{d}} {\bf N}}{\mathop{\mathrm{d}}
 \alpha}{\bf N}^{-1}\frac{\mathop{\mathrm{d}} {\bf N}}{\mathop{\mathrm{d}}
\alpha}{\bf N}^{-1} {\bf M}'- {\bf M}{\bf N}^{-1}\frac{\mathop{\mathrm{d}}^2
{\bf N}}{\mathop{\mathrm{d}} \alpha^2} {\bf N}^{-1}{\bf M}',
\end{align*}

\item
\begin{align*}
\frac{\mathop{\mathrm{d}}^2 \mbPsi^{+}(\alpha)}{\mathop{\mathrm{d}}
\alpha^2}\mbPsi(\alpha)\mbPsi^{+}(\alpha)  &=
\frac{\mathop{\mathrm{d}}^2 {\bf M}}{\mathop{\mathrm{d}}
\alpha^2}{\bf N}^{-1}{\bf M}' + {\bf M}{\bf N}^{-1}\frac{\mathop{\mathrm{d}}^2 {\bf
M}'}{\mathop{\mathrm{d}} \alpha^2} {\bf M}({\bf M}'{\bf M})^{-1}\\
&\quad + 2\frac{\mathop{\mathrm{d}} {\bf M}}{\mathop{\mathrm{d}} \alpha}{\bf N}^{-1}
\frac{\mathop{\mathrm{d}}{\bf M}'}{\mathop{\mathrm{d}} \alpha}{\bf M}({\bf M}'
{\bf M})^{-1}{\bf M}'\\
&\quad - 2\frac{\mathop{\mathrm{d}} {\bf M}}{\mathop{\mathrm{d}} \alpha}{\bf N}^{-1}
\frac{\mathop{\mathrm{d}} {\bf N}}{\mathop{\mathrm{d}} \alpha}{\bf N}^{-1}{\bf M}'\\
&\quad - 2{\bf M}{\bf N}^{-1}\frac{\mathop{\mathrm{d}} {\bf N}}{\mathop{\mathrm{d}}
\alpha}{\bf N}^{-1}
\frac{\mathop{\mathrm{d}} {\bf M}'}{\mathop{\mathrm{d}} \alpha}{\bf M}({\bf M}'
{\bf M})^{-1}{\bf M}',
\end{align*}

\item
\begin{align*}
\frac{\mathop{\mathrm{d}} \mbPsi^{+}(\alpha)}{\mathop{\mathrm{d}} \alpha}\mbPsi(\alpha)
\frac{\mathop{\mathrm{d}} \mbPsi^{+}(\alpha)}{\mathop{\mathrm{d}} \alpha}
&= \frac{\mathop{\mathrm{d}} {\bf M}}{\mathop{\mathrm{d}} \alpha}({\bf M}'{\bf M})^{-1}
{\bf M}'\frac{\mathop{\mathrm{d}} {\bf M}}{\mathop{\mathrm{d}} \alpha}{\bf N}^{-1}{\bf M}'
+ \frac{\mathop{\mathrm{d}} {\bf M}}{\mathop{\mathrm{d}} \alpha}{\bf N}^{-1}
\frac{\mathop{\mathrm{d}} {\bf M}'}{\mathop{\mathrm{d}} \alpha}\\
&\quad - \frac{\mathop{\mathrm{d}} {\bf M}}{\mathop{\mathrm{d}} \alpha}{\bf N}^{-1}
\frac{\mathop{\mathrm{d}} {\bf N}}{\mathop{\mathrm{d}} \alpha}{\bf N}^{-1}
\frac{\mathop{\mathrm{d}} {\bf M}'}{\mathop{\mathrm{d}} \alpha}\\
&\quad + {\bf M}{\bf N}^{-1}\frac{\mathop{\mathrm{d}} {\bf M}'}{\mathop{\mathrm{d}}
\alpha}{\bf M}({\bf M}'{\bf M})^{-1}{\bf N}({\bf M}'{\bf M})^{-1}{\bf M}'
\frac{\mathop{\mathrm{d}} {\bf M}}{\mathop{\mathrm{d}} \alpha}{\bf N}^{-1}{\bf M}'\\
&\quad + {\bf M}{\bf N}^{-1}\frac{\mathop{\mathrm{d}} {\bf M}'}{\mathop{\mathrm{d}}
\alpha}{\bf M}({\bf M}'{\bf M})^{-1}\frac{\mathop{\mathrm{d}} {\bf M}'}
{\mathop{\mathrm{d}} \alpha}\\ &\quad - {\bf M}{\bf N}^{-1}\frac{\mathop{\mathrm{d}}
{\bf M}'}{\mathop{\mathrm{d}} \alpha} {\bf M}({\bf M}'{\bf M})^{-1}
\frac{\mathop{\mathrm{d}} {\bf N}}{\mathop{\mathrm{d}} \alpha} {\bf N}^{-1}{\bf M}'\\
&\quad - {\bf M}{\bf N}^{-1}\frac{\mathop{\mathrm{d}} {\bf N}}{\mathop{\mathrm{d}}
 \alpha}({\bf M}'{\bf M})^{-1}{\bf M}'
\frac{\mathop{\mathrm{d}} {\bf M}}{\mathop{\mathrm{d}} \alpha}{\bf N}^{-1}{\bf M}'  \\
&\quad - {\bf M}{\bf N}^{-1}\frac{\mathop{\mathrm{d}} {\bf N}}{\mathop{\mathrm{d}} \alpha}
{\bf N}^{-1}\frac{\mathop{\mathrm{d}} {\bf M}'}{\mathop{\mathrm{d}} \alpha}
  +{\bf M}{\bf N}^{-1}\frac{\mathop{\mathrm{d}} {\bf N}}{\mathop{\mathrm{d}} \alpha}
{\bf N}^{-1}\frac{\mathop{\mathrm{d}} {\bf N}}{\mathop{\mathrm{d}} \alpha}{\bf N}^{-1}{\bf M}',
\end{align*}

\item
\begin{align*}
\frac{\mathop{\mathrm{d}} \mbPsi^{+}(\alpha)}{\mathop{\mathrm{d}} \alpha}
\frac{\mathop{\mathrm{d}} [{\bf M}({\bf M}'{\bf M})^{-1}{\bf M}']}{\mathop{\mathrm{d}}
\alpha}
&= \frac{\mathop{\mathrm{d}} {\bf M}}{\mathop{\mathrm{d}} \alpha}{\bf N}^{-1}{\bf M}'
\frac{\mathop{\mathrm{d}} {\bf M}}{\mathop{\mathrm{d}} \alpha}({\bf M}'{\bf M})^{-1}
{\bf M}' \\
&\quad - \frac{\mathop{\mathrm{d}} {\bf M}}{\mathop{\mathrm{d}} \alpha}{\bf N}^{-1}
 ( \frac{\mathop{\mathrm{d}} {\bf M}'}{\mathop{\mathrm{d}} \alpha}{\bf M}
 + {\bf M}'\frac{\mathop{\mathrm{d}} {\bf M}}{\mathop{\mathrm{d}} \alpha})({\bf M}'
{\bf M})^{-1}{\bf M}' \\
&\quad + \frac{\mathop{\mathrm{d}} {\bf M}}{\mathop{\mathrm{d}} \alpha}{\bf N}^{-1}
\frac{\mathop{\mathrm{d}} {\bf M}'}{\mathop{\mathrm{d}} \alpha}
+ {\bf M}{\bf N}^{-1}\frac{\mathop{\mathrm{d}} {\bf M}'}{\mathop{\mathrm{d}} \alpha}
\frac{\mathop{\mathrm{d}} {\bf M}}{\mathop{\mathrm{d}} \alpha}({\bf M}'{\bf M})^{-1}{\bf M}^{'}\\
&\quad \quad - {\bf M}{\bf N}^{-1} \frac{\mathop{\mathrm{d}} {\bf M}'}
{\mathop{\mathrm{d}} \alpha}\mbM(\mbM'\mbM)^{-1}
  (\frac{\mathop{\mathrm{d}} {\bf M}'}{\mathop{\mathrm{d}} \alpha}\mbM
+ \mbM' \frac{\mathop{\mathrm{d}} {\bf M}}{\mathop{\mathrm{d}} \alpha})(\mbM'\mbM)^{-1}\mbM'  \\
&\hspace{1cm} + \mbM\mbN^{-1} \frac{\mathop{\mathrm{d}} {\bf M}'}{\mathop{\mathrm{d}}\alpha}
\mbM(\mbM'\mbM)^{-1}\frac{\mathop{\mathrm{d}} {\bf M}'}{\mathop{\mathrm{d}} \alpha},
\end{align*}
\item
\begin{align*}
\frac{\mathop{\mathrm{d}} [{\bf M}({\bf M}'{\bf M})^{-1}{\bf M}']}{\mathop{\mathrm{d}} \alpha}
\frac{\mathop{\mathrm{d}} \mbPsi^{+}(\alpha)}{\mathop{\mathrm{d}} \alpha}
&= \frac{\mathop{\mathrm{d}} \mbM}{\mathop{\mathrm{d}} \alpha}(\mbM'\mbM)^{-1}\mbM'
\frac{\mathop{\mathrm{d}} \mbM}{\mathop{\mathrm{d}} \alpha}\mbN^{-1}\mbM'\\
&\quad - \mbM(\mbM'\mbM)^{-1} (\frac{\mathop{\mathrm{d}} \mbM'}{\mathop{\mathrm{d}} \alpha}\mbM
+ \mbM'\frac{\mathop{\mathrm{d}} \mbM}{\mathop{\mathrm{d}} \alpha})(\mbM'\mbM)^{-1}\mbM'
\frac{\mathop{\mathrm{d}} \mbM}{\mathop{\mathrm{d}} \alpha}\mbN^{-1}\mbM'\\
&\quad + \mbM(\mbM'\mbM)^{-1}\frac{\mathop{\mathrm{d}} \mbM'}{\mathop{\mathrm{d}} \alpha}
\frac{\mathop{\mathrm{d}} \mbM}{\mathop{\mathrm{d}} \alpha}\mbN^{-1}\mbM'
+\frac{\mathop{\mathrm{d}} \mbM}{\mathop{\mathrm{d}} \alpha}\mbN^{-1}
\frac{\mathop{\mathrm{d}} \mbM'}{\mathop{\mathrm{d}} \alpha} \\
\end{align*}
\begin{align*}\nonumber
&\quad - \mbM(\mbM'\mbM)^{-1} (\frac{\mathop{\mathrm{d}} \mbM'}{\mathop{\mathrm{d}} \alpha}\mbM
+ \mbM'\frac{\mathop{\mathrm{d}} \mbM}{\mathop{\mathrm{d}} \alpha})\mbN^{-1}
\frac{\mathop{\mathrm{d}} \mbM'}{\mathop{\mathrm{d}} \alpha}\\
&\quad + \mbM(\mbM'\mbM)^{-1}\frac{\mathop{\mathrm{d}} \mbM'}{\mathop{\mathrm{d}} \alpha}
\mbM\mbN^{-1}\frac{\mathop{\mathrm{d}} \mbM'}{\mathop{\mathrm{d}} \alpha}
- \frac{\mathop{\mathrm{d}} \mbM}{\mathop{\mathrm{d}} \alpha}\mbN^{-1}
\frac{\mathop{\mathrm{d}} \mbN}{\mathop{\mathrm{d}} \alpha}\mbN^{-1}\mbM'\\
&\quad + \mbM(\mbM'\mbM)^{-1}(\frac{\mathop{\mathrm{d}} \mbM'}{\mathop{\mathrm{d}} \alpha}
\mbM + \mbM'\frac{\mathop{\mathrm{d}} \mbM}{\mathop{\mathrm{d}} \alpha})\mbN^{-1}
\frac{\mathop{\mathrm{d}} \mbN}{\mathop{\mathrm{d}} \alpha}\mbN^{-1}\mbM'\\
&\quad - \mbM(\mbM'\mbM)^{-1}\frac{\mathop{\mathrm{d}} \mbM'}{\mathop{\mathrm{d}} \alpha}
\mbM\mbN^{-1}\frac{\mathop{\mathrm{d}} \mbN}{\mathop{\mathrm{d}} \alpha}\mbN^{-1}\mbM'.
\end{align*}
\end{enumerate}

\noindent Thus, by the results (A.4) and (i)-(vi) and some straightforward manipulations,
\begin{align}  \nonumber
  \mbM'\mbPsi^{+}(\alpha)\frac{\mathop{\mathrm{d}}^2 \mbPsi(\alpha)}{\mathop{\mathrm{d}}
  \alpha^2}\mbPsi^{+}(\alpha)\mbM
&=\mbM'\mbM \mbN^{-1}\frac{\mathop{\mathrm{d}}^2 \mbN}{\mathop{\mathrm{d}}
   \alpha^2}\mbN^{-1}\mbM'\mbM     \\  \nonumber
&\quad + 2\mbM'\frac{\mathop{\mathrm{d}} \mbM}{\mathop{\mathrm{d}} \alpha}(\mbM'\mbM)^{-1}\mbM'
 \frac{\mathop{\mathrm{d}} \mbM}{\mathop{\mathrm{d}} \alpha}\mbN^{-1}\mbM'\mbM   \\ \nonumber
&\quad + 2 \mbM'\mbM \mbN^{-1}\frac{\mathop{\mathrm{d}} \mbM'}
   {\mathop{\mathrm{d}} \alpha}\mbM(\mbM'\mbM)^{-1}
   \frac{\mathop{\mathrm{d}} \mbM'}{\mathop{\mathrm{d}} \alpha}\mbM    \\ \nonumber
&\quad + 2 \mbM'\mbM \mbN^{-1}\frac{\mathop{\mathrm{d}} \mbM'}{\mathop{\mathrm{d}} \alpha}
  \mbM(\mbM'\mbM)^{-1}\mbN (\mbM'\mbM)^{-1}\mbM'\frac{\mathop{\mathrm{d}} \mbM}
  {\mathop{\mathrm{d}} \alpha}\mbN^{-1} \mbM'\mbM    \\ \nonumber
&\quad + 2 \mbM'\mbM \mbN^{-1}\frac{\mathop{\mathrm{d}}\mbM'}{\mathop{\mathrm{d}} \alpha}
  \mbM(\mbM'\mbM)^{-1}\mbM'\frac{\mathop{\mathrm{d}} \mbM}{\mathop{\mathrm{d}} \alpha}  \\  \nonumber
&\quad + 2\frac{\mathop{\mathrm{d}} \mbM'}{\mathop{\mathrm{d}} \alpha}\mbM(\mbM'\mbM)^{-1}
   \mbM^{-1}\frac{\mathop{\mathrm{d}} \mbM}{\mathop{\mathrm{d}}\alpha}\mbN^{-1} \mbM'\mbM \\ \nonumber
&\quad - \mbM'\frac{\mathop{\mathrm{d}}^2 \mbM}{\mathop{\mathrm{d}} \alpha^2}\mbN^{-1}
   \mbM'\mbM -\mbM'\mbM \mbN^{-1}\frac{\mathop{\mathrm{d}}^2 \mbM'}{\mathop{\mathrm{d}}
   \alpha^2}\mbM   \\ \nonumber
&\quad - 2\frac{\mathop{\mathrm{d}} \mbM'}{\mathop{\mathrm{d}} \alpha}
   \frac{\mathop{\mathrm{d}} \mbM}{\mathop{\mathrm{d}} \alpha}\mbN^{-1} \mbM'\mbM
  - 2 \mbM'\mbM \mbN^{-1}\frac{\mathop{\mathrm{d}} \mbM'}{\mathop{\mathrm{d}} \alpha}
   \frac{\mathop{\mathrm{d}} \mbM}{\mathop{\mathrm{d}} \alpha}  \\ \nonumber
&\quad - 2 \mbM'\mbM\mbN^{-1}\frac{\mathop{\mathrm{d}} \mbM'}{\mathop{\mathrm{d}}
   \alpha}\mbM(\mbM'\mbM)^{-1}\frac{\mathop{\mathrm{d}} \mbN}{\mathop{\mathrm{d}}
   \alpha}\mbN^{-1} \mbM'\mbM  \\
&\quad - 2 \mbM'\mbM \mbN^{-1}\frac{\mathop{\mathrm{d}} \mbN}{\mathop{\mathrm{d}}
   \alpha}(\mbM'\mbM)^{-1}\mbM'\frac{\mathop{\mathrm{d}} \mbM}{\mathop{\mathrm{d}}
   \alpha}\mbN^{-1} \mbM'\mbM.
\end{align}

\noindent Also, note that
\[
\frac{\mathop{\mathrm{d}} \mbPsi^{+}(\alpha)}{\mathop{\mathrm{d}} \alpha}
= \frac{\mathop{\mathrm{d}} \mbM}{\mathop{\mathrm{d}} \alpha}\mbN^{-1}\mbM'
+ \mbM\mbN^{-1}\frac{\mathop{\mathrm{d}} \mbM'}{\mathop{\mathrm{d}} \alpha}
- \mbM\mbN^{-1}\frac{\mathop{\mathrm{d}} \mbN}{\mathop{\mathrm{d}} \alpha}\mbN^{-1}\mbM'.
\]
Thus, the stationary point of $\mbPsi^{+}(\alpha)$ satisfies the following equation
\[
\frac{\mathop{\mathrm{d}} \mbM}{\mathop{\mathrm{d}} \alpha}\mbN^{-1}\mbM'
+ \mbM\mbN^{-1}\frac{\mathop{\mathrm{d}} \mbM'}{\mathop{\mathrm{d}} \alpha}
= \mbM\mbN^{-1}\frac{\mathop{\mathrm{d}} \mbN}{\mathop{\mathrm{d}} \alpha}\mbN^{-1}\mbM',
\]
which implies that
\[
(\mbM'\mbM)^{-1}\mbM\frac{\mathop{\mathrm{d}} \mbM}{\mathop{\mathrm{d}} \alpha}\mbN^{-1}
+ \mbN^{-1}\frac{\mathop{\mathrm{d}} \mbM'}{\mathop{\mathrm{d}} \alpha}
\mbM(\mbM'\mbM)^{-1}
= \mbN^{-1}\frac{\mathop{\mathrm{d}} \mbN}{\mathop{\mathrm{d}} \alpha}\mbN^{-1}.
\]
Namely,
\[
 \mbM'\mbM\mbN^{-1}\frac{\mathop{\mathrm{d}} \mbN}{\mathop{\mathrm{d}} \alpha}
  - \mbM'\frac{\mathop{\mathrm{d}} \mbM}{\mathop{\mathrm{d}} \alpha}
  =\mbM'\mbM \mbN^{-1}\frac{\mathop{\mathrm{d}} \mbM'}{\mathop{\mathrm{d}} \alpha}
\mbM(\mbM'\mbM)^{-1}\mbN.
\]
Thus,
\begin{align*}
& (\mbM'\mbM\mbN^{-1}\frac{\mathop{\mathrm{d}} \mbN}{\mathop{\mathrm{d}} \alpha}
- \mbM'\frac{\mathop{\mathrm{d}} \mbM}{\mathop{\mathrm{d}} \alpha}) \mbN^{-1}
 (\mbM'\mbM \mbN^{-1}\frac{\mathop{\mathrm{d}} \mbN}{\mathop{\mathrm{d}} \alpha}
- \mbM'\frac{\mathop{\mathrm{d}} \mbM}{\mathop{\mathrm{d}} \alpha})'\\
&\qquad =\mbM'\mbM \mbN^{-1}\frac{\mathop{\mathrm{d}} \mbM'}{\mathop{\mathrm{d}} \alpha}
  \mbM(\mbM'\mbM)^{-1}\mbN(\mbM'\mbM)^{-1}\mbM'\frac{\mathop{\mathrm{d}} \mbM}
  {\mathop{\mathrm{d}}\alpha}\mbN^{-1} \mbM'\mbM.
\end{align*}
Furthermore, ${\bf B}^{+}(\mbS)$ is convex on $\mbS$,
thus $\frac{\mathop{\mathrm{d}}^2 \mbPsi^{+}(\alpha)}{\mathop{\mathrm{d}} \alpha^2}
\succeq {\bf 0}$, that is,
\begin{align*}
&\frac{\mathop{\mathrm{d}}^2 \mbM}{\mathop{\mathrm{d}} \alpha^2}\mbN^{-1}\mbM'
+ \mbM\mbN^{-1}\frac{\mathop{\mathrm{d}}^2 \mbM'}{\mathop{\mathrm{d}} \alpha^2}
- \mbM\mbN^{-1}\frac{\mathop{\mathrm{d}}^2 \mbN}
{\mathop{\mathrm{d}} \alpha^2}\mbN^{-1}\mbM'\\
&\qquad \succeq 2( \mbM\mbN^{-1}\frac{\mathop{\mathrm{d}} \mbN}
{\mathop{\mathrm{d}} \alpha}
 - \frac{\mathop{\mathrm{d}} \mbM}{\mathop{\mathrm{d}} \alpha}) \mbN^{-1}
 ( \mbM\mbN^{-1}\frac{\mathop{\mathrm{d}} \mbN}{\mathop{\mathrm{d}} \alpha}
  - \frac{\mathop{\mathrm{d}} \mbM}{\mathop{\mathrm{d}} \alpha})'.
\end{align*}
Substitute these results into (A.5), then
\begin{align*}
&\left. \mbM'\mbPsi^{+}(\alpha)\frac{\mathop{\mathrm{d}}^2 \mbPsi(\alpha)}
{\mathop{\mathrm{d}}\alpha^2}\mbPsi^{+}(\alpha)\mbM \right|_{\alpha = \alpha_0}\\
&\qquad \preceq 2 [\frac{\mathop{\mathrm{d}} \mbM'}{\mathop{\mathrm{d}} \alpha}
\mbM(\mbM'\mbM)^{-1}\mbM'\frac{\mathop{\mathrm{d}} \mbM}{\mathop{\mathrm{d}}
\alpha}\mbN^{-1} \mbM'\mbM + \mbM'\mbM \mbN^{-1}\frac{\mathop{\mathrm{d}} \mbM'}
{\mathop{\mathrm{d}}\alpha}\mbM(\mbM'\mbM)^{-1}\mbM'
\frac{\mathop{\mathrm{d}} \mbM}{\mathop{\mathrm{d}} \alpha} ]\\
&\qquad\qquad -2 [ \frac{\mathop{\mathrm{d}} \mbM'}{\mathop{\mathrm{d}} \alpha}
\frac{\mathop{\mathrm{d}} \mbM}{\mathop{\mathrm{d}} \alpha}\mbN^{-1}\mbM'\mbM
+ \mbM'\mbM \mbN^{-1}\frac{\mathop{\mathrm{d}} \mbM'}{\mathop{\mathrm{d}} \alpha}
\frac{\mathop{\mathrm{d}} \mbM}{\mathop{\mathrm{d}} \alpha} ]\\
&\qquad \preceq {\bf 0}.
\end{align*}
Thus,
\[
\left. \frac{\mathop{\mathrm{d}}^2 \mbPsi(\alpha)}{\mathop{\mathrm{d}} \alpha^2}
\right|_{\alpha = \alpha_0} \preceq {\bf 0}.
\]
By Lemma 3, then
\[
\frac{\mathop{\mathrm{d}}^2 \mbPsi(\alpha)}{\mathop{\mathrm{d}} \alpha^2} \preceq
\left. \frac{\mathop{\mathrm{d}}^2 \mbPsi(\alpha)}{\mathop{\mathrm{d}} \alpha^2}
\right|_{\alpha = \alpha_0} \preceq {\bf 0}, ~\forall ~\alpha \in (0,1).
\]
Therefore, ${\bf B}(\mbS)$ is concave on ${\cal S}$.

\vspace{0.3cm}
\noindent {\bf Lemma 5}.
{\it Let ${\bf A}_i$, $i = 1$, $2$, be $p \times p$ p.s.d. of rank $r$ ($r \leq p$).
Also let ${\bf D}_r = \mathop{\mathrm{diag}} (d_1, \cdots, d_r)$
with elements being the non-zero eigenvalues of ${\bf A}_2{\bf A}_1^{+}$,
where ${\bf A}_1^{+}$ denotes the Moore-Penrose generalized inverse of ${\bf A}_1$.
Then there exists a nonsingular matrix $G$ such that
\[
{\bf A}_1 = {\bf G} \left[
\begin{array}{cc}
{\bf I}_r& {\bf 0}\\
{\bf 0}& {\bf 0}
\end{array}
\right]
{\bf G}' \quad \mbox{and} \quad
{\bf A}_2 = {\bf G} \left[
\begin{array}{cc}
{\bf D}_r& {\bf 0}\\
{\bf 0}& {\bf 0}
\end{array}
\right]
{\bf G}'
\]}

\indent {\bf Proof.} By Theorem A.4.1 of Anderson [1],
\[
{\bf A}_1 = {\bf F} \left[
\begin{array}{cc}
{\bf I}_r& {\bf 0}\\
{\bf 0}& {\bf 0}
\end{array}
\right]
{\bf F}',
\]
where ${\bf F}$ is a nonsingular matrix.
Let ${\bf A}^{\star} = {\bf F}^{-1}{\bf A}_2({\bf F}^{-1})'$, then
${\bf A}^{\star}$ is a p.s.d. with rank $r$.
Write
\[
{\bf A}^{\star}
= \left(
\begin{array}{cc}
{\bf A}_{11}^{\star}& {\bf A}_{12}^{\star}\\
{\bf A}_{21}^{\star}& {\bf A}_{22}^{\star}
\end{array}
\right),
\]
and take
\[
{\bf C} = \left[
\begin{array}{cc}
{\bf C}_{11}& {\bf 0}\\
{\bf 0}& {\bf C}_{22}
\end{array}
\right]_,
\]
where ${\bf C}_{11} \in \mathcal{Q}(r)$, ${\bf C}_{22} \in \mathcal{Q}(p-r)$
such that ${\bf C}_{22}'{\bf A}_{22}^{\star} = 0$.
Thus, there exists a matrix ${\bf C} \in \mathcal{Q}(p)$,
the group of $p \times p$ orthogonal matrices such that
\[
{\bf C}'{\bf A}^{\star}{\bf C}
= \left[
\begin{array}{cc}
{\bf D}_r& {\bf 0}\\
{\bf 0}& {\bf 0}
\end{array}
\right]
\]
that is,
\[
\left\{
\begin{array}{lll}
{\bf C}'{\bf F}^{-1}{\bf A}_1({\bf F}^{-1})'{\bf C}
&=& \left[
\begin{array}{cc}
{\bf D}_r& {\bf 0}\\
{\bf 0}& {\bf 0}
\end{array}
\right]\\
{\bf C}'{\bf F}^{-1}{\bf A}_2({\bf F}^{-1})'{\bf C}
&=& \left[
\begin{array}{cc}
{\bf I}_r& {\bf 0}\\
{\bf 0}& {\bf 0}
\end{array}
\right]
\end{array}
\right._.
\]
Let ${\bf G} = {\bf FC}$, thus
\[
{\bf A}_1 = {\bf G} \left[
\begin{array}{cc}
{\bf I}_r& {\bf 0}\\
{\bf 0}& {\bf 0}
\end{array}
\right]
{\bf G}'
\]
and
\[
{\bf A}_2 = {\bf G} \left[
\begin{array}{cc}
{\bf D}_r& {\bf 0}\\
{\bf 0}& {\bf 0}
\end{array}
\right]
{\bf G}'.
\]

\vspace{0.3cm}
\indent {\bf Lemma 6}.
{\it Let ${\bf A}$ be an $p \times p$ p.s.d. matrix and ${\bf x}$ be a $p \times 1$ vector.
Let ${\bf A} = {\bf BCB}'$ and denotes ${\bf A}^{-} = ({\bf B}')^{+}{\bf C}^{+}{\bf B}^{+}$,
where ${\bf D}^{+}$ denotes the Moore-Penrose generalized inverse of ${\bf D}$.
If ${\bf A}^{-}$ is a generalized inverse of ${\bf A}$, then
$f({\bf x}, {\bf A}) = {\bf x}'{\bf A}^{-}{\bf x}$ is convex on $R^{p}\times {\cal S}$}.

\indent {\bf Proof.} Since $f$ is continuous in $({\bf x}, {\bf A})$, it suffices to show that
\[
({\bf x} + {\bf y})'({\bf A}_1 + {\bf A}_2)^{-}({\bf x} + {\bf y})
\leq {\bf x}'{\bf A}_1^{-}{\bf x} + {\bf y}'{\bf A}_2^{-}{\bf y}.
\]
Write ${\bf A} = {\bf A}^{1/2}({\bf A}^{1/2})'$ and take
\begin{align*}
{\bf u} &= ({\bf A}_1^\frac{1}{2})^{-}{\bf x} - ({\bf A}_1^\frac{1}{2})'({\bf A}_1
+ {\bf A}_2)^{-}({\bf x} + {\bf y})\\
{\bf v} &= ({\bf A}_2^\frac{1}{2})^{-}{\bf y} - ({\bf A}_2^\frac{1}{2})'({\bf A}_1
+ {\bf A}_2)^{-}({\bf x} + {\bf y}).
\end{align*}
By Lemma 5, ${\bf A}_1 = {\bf G} \left[ \begin{array}{cc} {\bf I}_r& {\bf 0}\\
{\bf 0}& {\bf 0} \end{array} \right] {\bf G}'$
and ${\bf A}_2 = {\bf G} \left[ \begin{array}{cc} {\bf D}_r& {\bf 0}\\
{\bf 0}& {\bf 0} \end{array} \right] {\bf G}'$,
where ${\bf G}$ is nonsingular and ${\bf D}_r = \mathop{\mathrm{diag}}(d_1, \cdots, d_r)$
with $d_i$ being the non-zero eigenvalues of ${\bf A}_2{\bf A}_1^{+}$.
Write ${\bf A}_1^\frac{1}{2} = {\bf G} \left[ \begin{array}{cc} {\bf I}_r
&{\bf 0}\\ {\bf 0}& {\bf 0} \end{array} \right]$
and ${\bf A}_2^\frac{1}{2} = {\bf G} \left[ \begin{array}{cc} {\bf D}_r& {\bf 0}\\
{\bf 0}& {\bf 0} \end{array} \right]$, then
\[
({\bf A}_1^\frac{1}{2})^{-} = \left[
\begin{array}{cc}
{\bf I}_r& {\bf 0}\\
{\bf 0}& {\bf 0}
\end{array}
\right]
{\bf G}^{-1} \quad \mbox{and} \quad
({\bf A}_2^\frac{1}{2})^{-} = \left[
\begin{array}{cc}
{\bf D}_r^\frac{-1}{2}& {\bf 0}\\
{\bf 0}& {\bf 0}
\end{array}
\right]
{\bf G}^{-1}~_.
\]
Notably,
\begin{align*}
0 &< {\bf u}'{\bf u} + {\bf v}'{\bf v}\\
  &= {\bf x}'\big( ({\bf A}_1^\frac{1}{2})^{-}\big)'({\bf A}_1^\frac{1}{2})^{-}{\bf x}
    -{\bf x}'\big( ({\bf A}_1^\frac{1}{2})^{-} \big)'({\bf A}_1^\frac{1}{2})'({\bf A}_1
   + {\bf A}_2)^{-}({\bf x} + {\bf y})\\
  & -({\bf x} + {\bf y})'({\bf A}_1 + {\bf A}_2)^{-}{\bf A}_1^\frac{1}{2}
  ({\bf A}_1^\frac{1}{2})^{-}{\bf x}+({\bf x} + {\bf y})'({\bf A}_1 + {\bf A}_2)^{-}{\bf
    A}_1^\frac{1}{2}({\bf A}_1^\frac{1}{2})'({\bf A}_1 + {\bf A}_2)^{-}({\bf x}
   + {\bf y})\\
  & + {\bf y}' \big( ({\bf A}_2^\frac{1}{2})^{-} \big)'({\bf A}_2^\frac{1}{2})^{-}{\bf y}
    - {\bf y}' \big( ({\bf A}_2^\frac{1}{2})^{-} \big)'({\bf A}_2^\frac{1}{2})'
      ({\bf A}_1 + {\bf A}_2)^{-}({\bf x} + {\bf y})\\
  & -({\bf x} + {\bf y})'({\bf A}_1 + {\bf A}_2)^{-}({\bf A}_2^\frac{1}{2})
   ({\bf A}_2^\frac{1}{2})^{-}{\bf y} + ({\bf x} + {\bf y})'({\bf A}_1 + {\bf A}_2)^{-}
   \big( {\bf A}_2^\frac{1}{2}({\bf A}_2^\frac{1}{2})'\big)
      ({\bf A}_1 + {\bf A}_2)^{-}({\bf x} + {\bf y}).
\end{align*}
Moreover,
\begin{align*}
\big( ({\bf A_1}^\frac{1}{2})^{-} \big)'({\bf A}_1^\frac{1}{2})'({\bf A}_1
+ {\bf A}_2)^{-}
&= ({\bf G}^{-1})' \left[
\begin{array}{cc}
{\bf I}_r& {\bf 0}\\
{\bf 0}& {\bf 0}
\end{array}
\right]
\left[
\begin{array}{cc}
{\bf I}_r& {\bf 0}\\
{\bf 0}& {\bf 0}
\end{array}
\right] {\bf G}'({\bf G}^{-1})' \left[
\begin{array}{cc}
({\bf I}_r + {\bf D}_r)^{-1}& {\bf 0}\\
{\bf 0}& {\bf 0}
\end{array}
\right]
{\bf G}^{-1}\\
&= ({\bf G}^{-1})' \left[
\begin{array}{cc}
({\bf I}_r + {\bf D}_r)^{-1}& {\bf 0}\\
{\bf 0}& {\bf 0}
\end{array}
\right]
{\bf G}^{-1}\\
&= ({\bf A}_1 + {\bf A}_2)^{-}
\end{align*}
and
\[
\big( ({\bf A}^\frac{1}{2})^{-} \big)'({\bf A}_1^\frac{1}{2})^{-} = ({\bf A}_1^{-})'
= {\bf A}_1^{-}.
\]
Thus,
\[
0 < {\bf u}'{\bf u} + {\bf v}'{\bf v} = {\bf x}'{\bf A}_1^{-}{\bf x}
+ {\bf y}'{\bf A}_2^{-}{\bf y} + ({\bf x} + {\bf y})'
({\bf A}_1 + {\bf A}_2)^{-}({\bf x} + {\bf y}),
\]
and hence the lemma follows.

\vspace{0.5cm}
\noindent {\bf Acknowledgments}
\vspace{0.3cm}

\indent The author is grateful to Professor C. R. Rao for his helpful comments.
The work was partly supported by Grants from National Science Council of the Republic
of China under Contract No. NSC 93-2118-M-001-027.

\vspace {0.5cm}
\noindent {\bf  References}
\begin{enumerate}

\item T.W. Anderson, {An Introduction to Multivariate Statistical Analysis},
2nd edition. New York: Wiley, 1984.

\item A. Birnbaum, Characterization of complete classes of tests of some
multiparametric hypotheses, with application to likelihood ratio tests, {Ann. Math.
Statist.} {26} (1955) 21-36.

\item M.L. Eaton, A complete class theorem for multidimensional one sided
alternatives, {Ann. Math. Statist.} {41} (1970) 1884-1888.


\item J. Kiefer, R. Schwartz, Admissible Bayes character of $T^2$-, $R^2$-, and
other fully invariant tests for classical multivariate normal problems, {Ann. Math.
Statist.} {36} (1965) 747-770.

\item E.L. Lehmann, {Testing Statistical Hypotheses}, 2nd edition. New York:
Wiley, 1986.

\item J. Marden, M.D. Perlman, Invariant tests for means with covariates, {
Ann. Statist.} {8} (1980) 25-63.

\item A.W. Marshall, I. Olkin, {Inequalities: Theory of Majorization and Its
Applications}, New York: Academic Press, 1979.


\item J. Oosterhoff,  {Combination of One-sided Statistical Tests}, Amsterdam:
Mathematisch Centrum, 1969.


\item C.R. Rao, Tests with discriminant functions in multivariate analysis,
{Sankhy$\bar{a}$} {7} (1946) 407-413.

\item C.R. Rao, On some problems arising out of discrimination with multiple characters,
{Sankhy$\bar{a}$} {9} (1949) 343-366.

\item C.R. Rao, S.K. Mitra, {Generalized Inverse of Matrices and its Applications},
New York: Wiley, 1971.

\item S.N. Roy, On a heuristic method of test construction and its use in multivariate
analysis, {Ann. Math. Statist.} {24} (1953) 220-238.

\item R. Schwartz, Admissible tests in multivariate analysis of variance, {Ann. Math. Statist.}
{38} (1967) 698-710.


\item J.B. Simaika, On an optimum property of two important statistical tests,
{Biometrika}  {32}  (1941) 70-80.


\item C. Stein, The admissibility of Hotelling's $T^{2}$-test, {Ann. Math.
Statist.} {27} (1956) 616-623.


\end{enumerate}


\end{document}